\documentclass[12pt]{amsart}

 \usepackage{amsmath,amsthm,amsfonts,amssymb,verbatim}
 \usepackage{url}
 \usepackage{graphicx}
 \usepackage[all]{xy}
 \usepackage{wrapfig}
 \usepackage{picins}

\newcommand{\sC}{\mathcal{C}}

\newcommand{\bN}{\mathbb{N}}
\newcommand{\bZ}{\mathbb{Z}}
\newcommand{\bQ}{\mathbb{Q}}

\newcommand{\al}{\alpha}
\newcommand{\be}{\beta}

\newcommand{\Del}{\Delta}
\newcommand{\si}{\sigma}

\newcommand{\ten}{\otimes}
\newcommand{\pls}{\oplus}

\newcommand{\vm}{v_-}
\newcommand{\vp}{v_+}
\newcommand{\vmm}{\vm\ten\vm}
\newcommand{\vmp}{\vm\ten\vp}
\newcommand{\vpm}{\vp\ten\vm}

\newcommand{\longto}{\longrightarrow}
\newcommand{\longmap}{\longmapsto}

\newcommand{\unknot}
{\raisebox{-0.055in}{\includegraphics[scale=0.125]{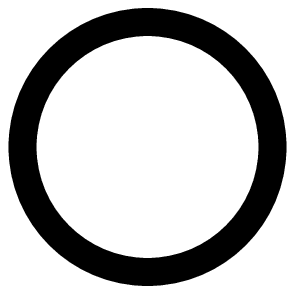}}}
\newcommand{\positive}
{\raisebox{-0.055in}
{\includegraphics[scale=0.125]{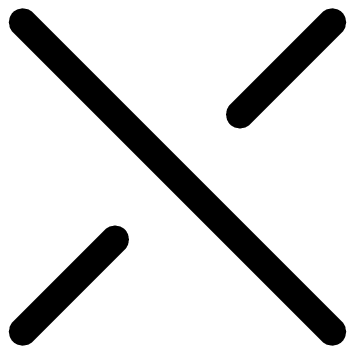}}}
\newcommand{\negative}
{\raisebox{-0.055in}
{\includegraphics[scale=0.125]{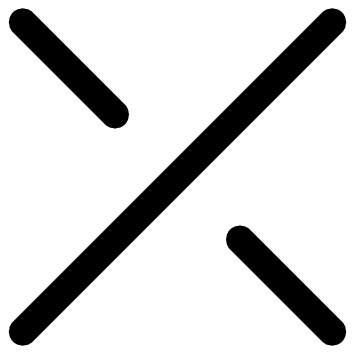}}}
\newcommand{\Cpositive}
{\raisebox{-0.08in}
{\includegraphics[scale=0.125]{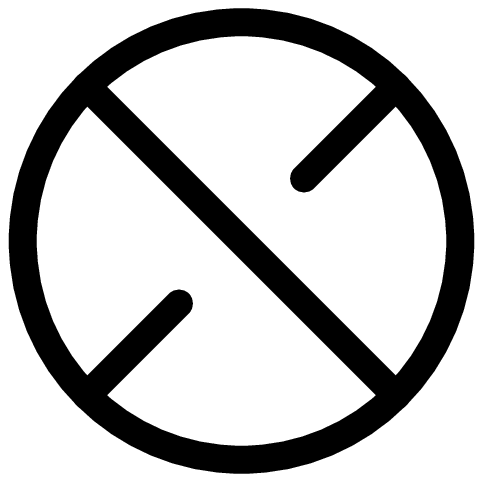}}}
\newcommand{\Cnegative}
{\raisebox{-0.08in}
{\includegraphics[scale=0.125]{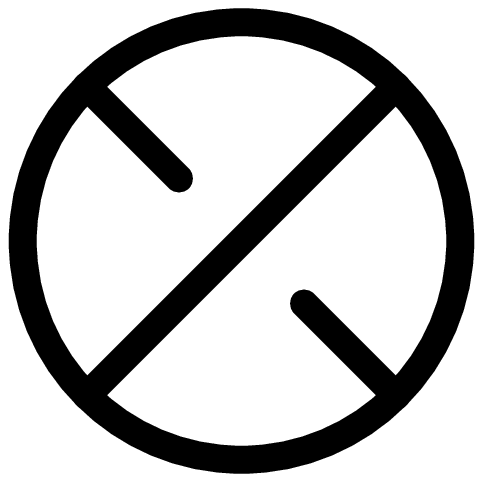}}}

\newcommand{\otherleftcross}
{\raisebox{-0.055in}
{\rotatebox{90}
{\includegraphics[scale=0.125]{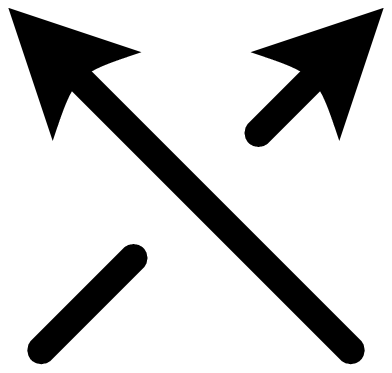}}}}
\newcommand{\otherrightcross}
{\raisebox{-0.055in}
{\rotatebox{90}
{\includegraphics[scale=0.125]{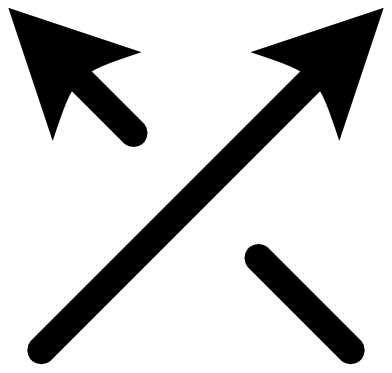}}}}
\newcommand{\zero}
{\raisebox{-0.055in}
{\includegraphics[scale=0.125]{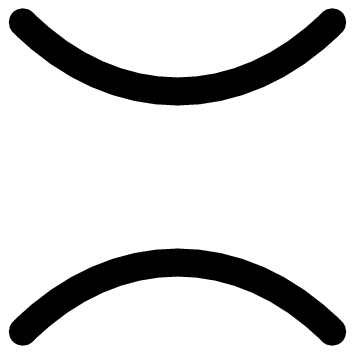}}}
\newcommand{\one}
{\raisebox{-0.055in}
{\includegraphics[scale=0.125]{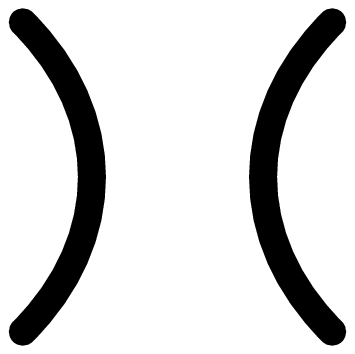}}}
\newcommand{\Czero}
{\raisebox{-0.08in}
{\includegraphics[scale=0.125]{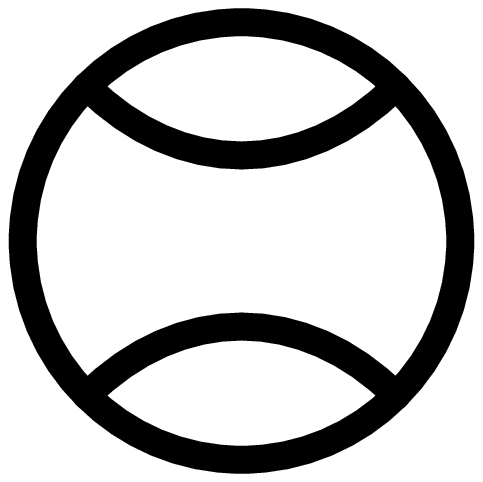}}}
\newcommand{\Cone}
{\raisebox{-0.08in}
{\includegraphics[scale=0.125]{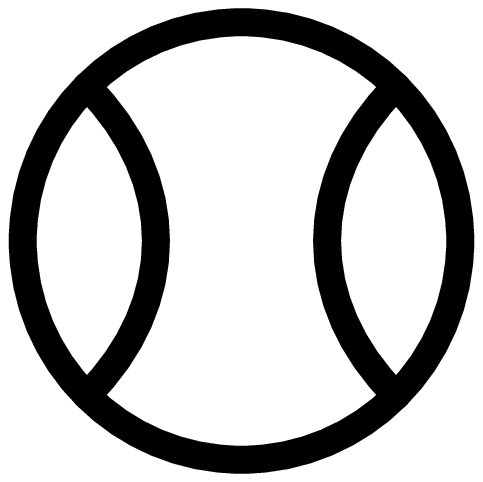}}}
\newcommand{\Gone}
{\raisebox{-0.08in}
{\includegraphics[scale=0.125]{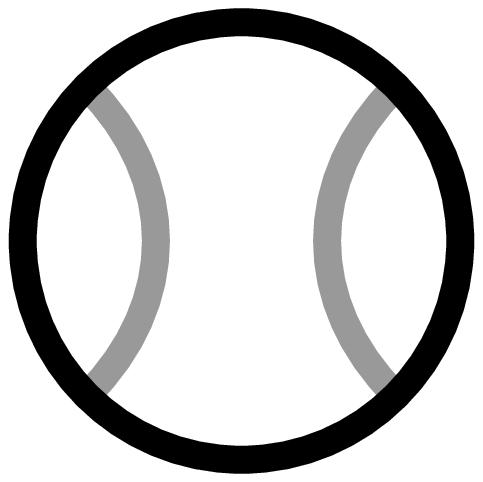}}}
\newcommand{\Gzero}
{\raisebox{-0.08in}
{\includegraphics[scale=0.125]{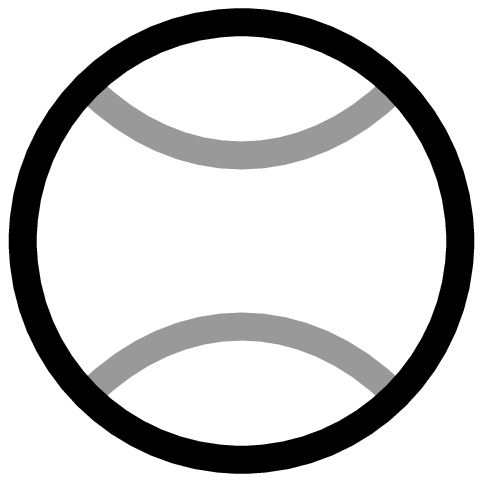}}}
\newcommand{\Gcross}
{\raisebox{-0.08in}
{\includegraphics[scale=0.125]{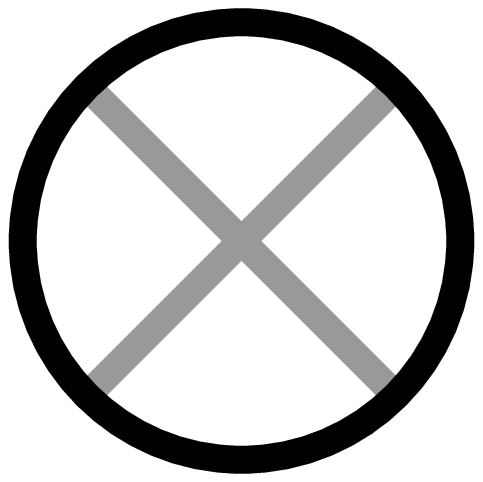}}}

\newcommand{\sigright}
{\raisebox{-0.08in}
{\includegraphics[scale=0.125]{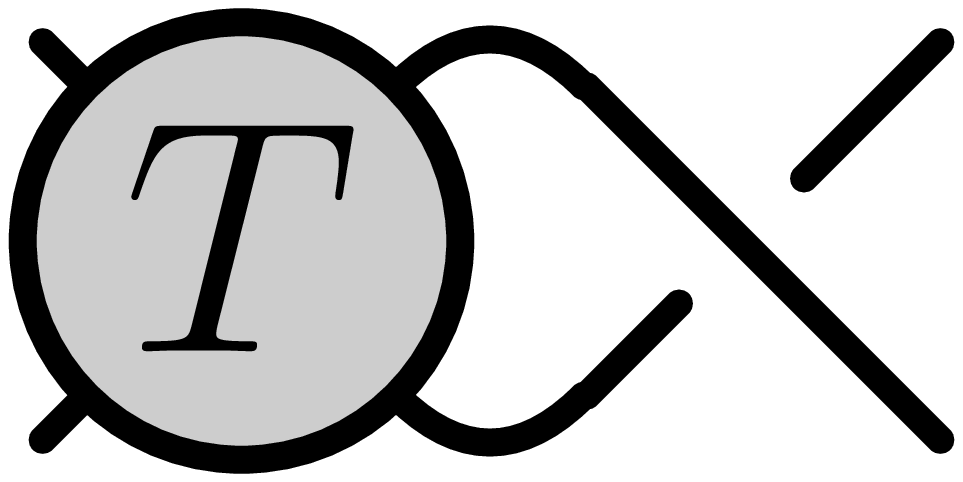}}}
\newcommand{\sigleft}
{\raisebox{-0.08in}
{\includegraphics[scale=0.125]{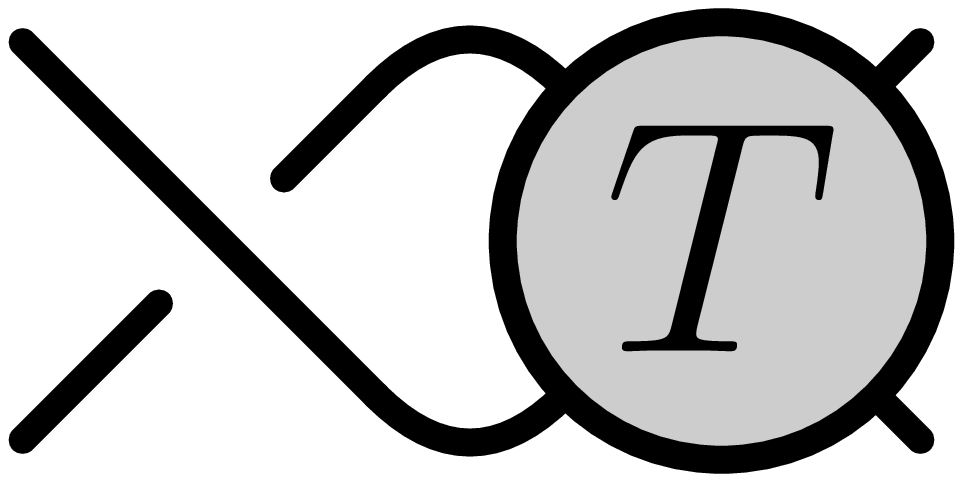}}}
\newcommand{\SIGright}
{\raisebox{-0.08in}
{\includegraphics[scale=0.125]{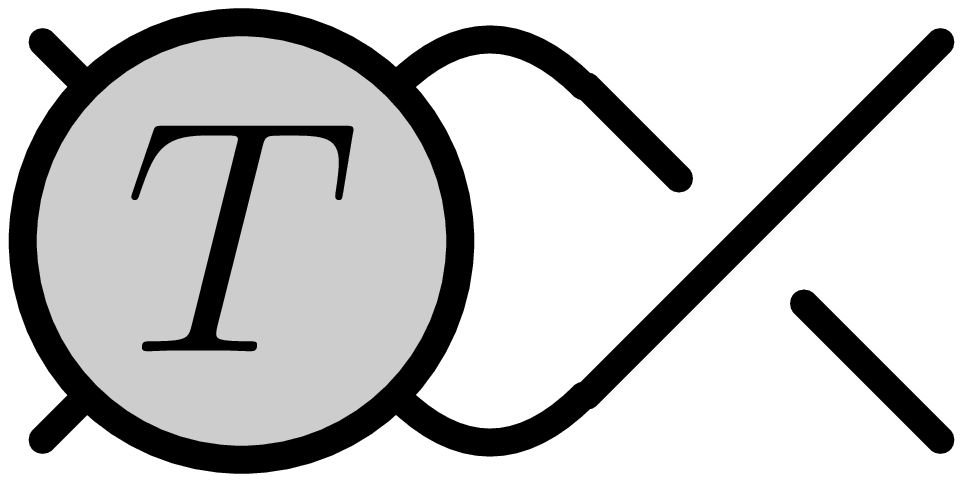}}}
\newcommand{\SIGleft}
{\raisebox{-0.08in}
{\includegraphics[scale=0.125]{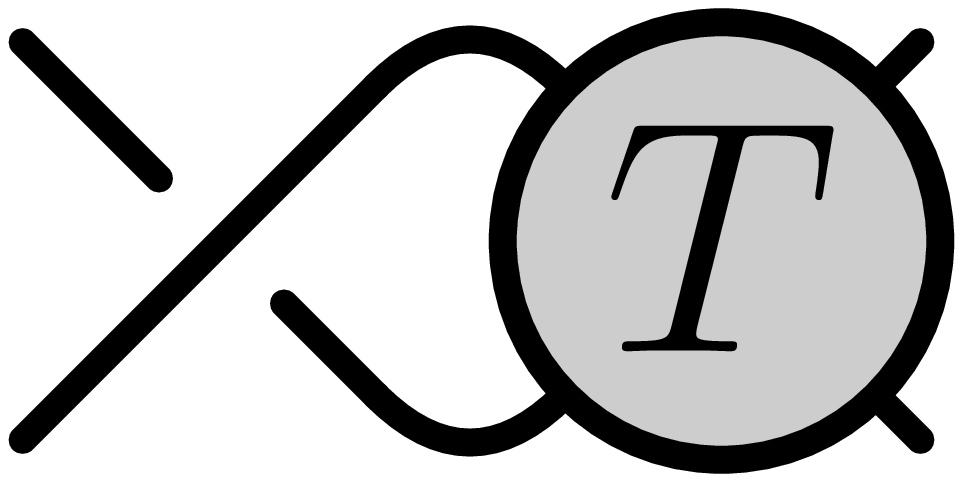}}}

\newtheorem{theorem}{Theorem}[section]

\newtheorem{remark}[theorem]{Remark}
\newtheorem{lemma}[theorem]{Lemma}

\newtheorem*{namedtheorem}{\theoremname}
\newcommand{\theoremname}{testing}

\title{Knots with identical Khovanov homology}
\author{Liam Watson}
\thanks{Supported by a NSERC Canadian Graduate Scholarship.}
\address{D\'epartement de Math\'ematiques, Universit\'e du Qu\'ebec \`a Montr\'eal,  Montr\'eal Canada, H3C 3P8}
\email{liam@math.uqam.ca}
\date{January 9, 2007}

\begin{document}

\begin{abstract}
We give a recipe for constructing families of distinct knots that have identical Khovanov homology and give examples of pairs of prime knots, as well as infinite families, with this property.
\end{abstract}

\maketitle

\section{Introduction}\label{Introduction}

The aim of this note is to present a construction giving rise to distinct knots that cannot be distinguished using Khovanov homology. Our main tool is the long exact sequence in Khovanov homology which is presented, along with a review of Khovanov homology, in section \ref{Notation}. In section \ref{Particular Example} we present a detailed calculation for a particular pair of knots, and in section \ref{Construction} we present a general construction for producing knots that can be handled via a similar calculation. This calculation, in the general setting, is carried out in section \ref{Proof}. In section \ref{Examples} we show how to apply this construction to obtain pairs of distinct prime knots with identical Khovanov homology (theorem \ref{thm:pairs}). These examples are distinguished by the HOMFLYPT polynomial, and as such must have distinct triply-graded link homology \cite{Khovanov2005}. We also give a construction of infinite families of distinct knots with identical Khovanov homology (theorem \ref{thm:infinite}). Finally, in section \ref{Mutants} we conclude with some remarks on a particular family of mutants.

\section{Notation}\label{Notation}

\parpic[r]{$\begin{array}{c}\xymatrix@R=0pt@C=5pt{
 && {\begin{array}{c}\includegraphics[scale=0.09375]{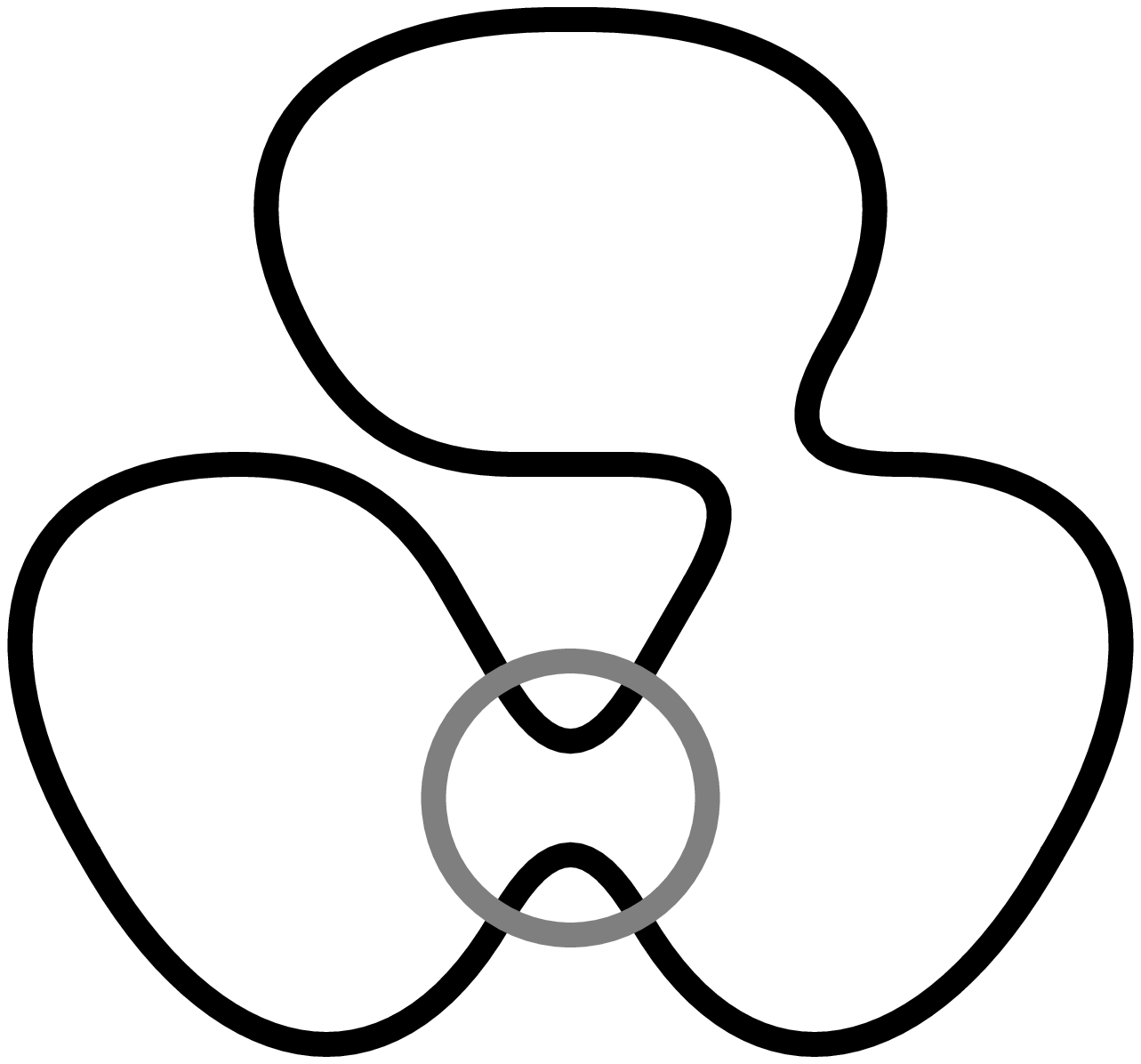}\\[-5pt] \textrm{\tiny{(1,0,0)}} \end{array}}\ar@{..}[d]\ar[rr]^-{-\Delta}\ar@{-}[dr] && {\begin{array}{c}\includegraphics[scale=0.09375]{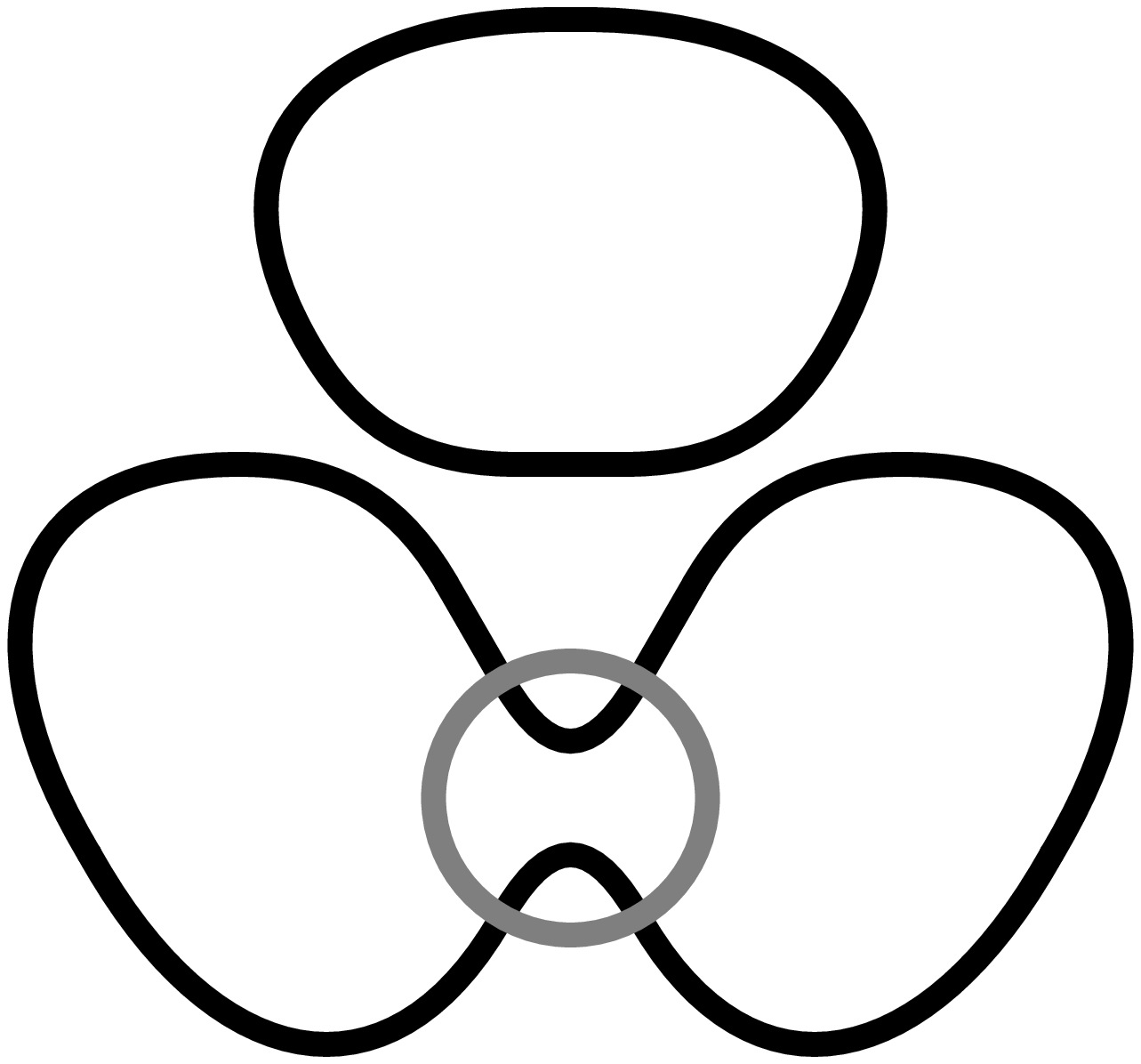}\\[-5pt] \textrm{\tiny{(1,1,0)}} \end{array}}\ar@{..}[d]\ar[ddrr]^-{\Delta} && \\
 && {\oplus}\ar@{..}[d] & \ar[dr]_(0.3){-\Delta} & {\oplus}\ar@{..}[d] && \\
 {\begin{array}{c}\includegraphics[scale=0.09375]{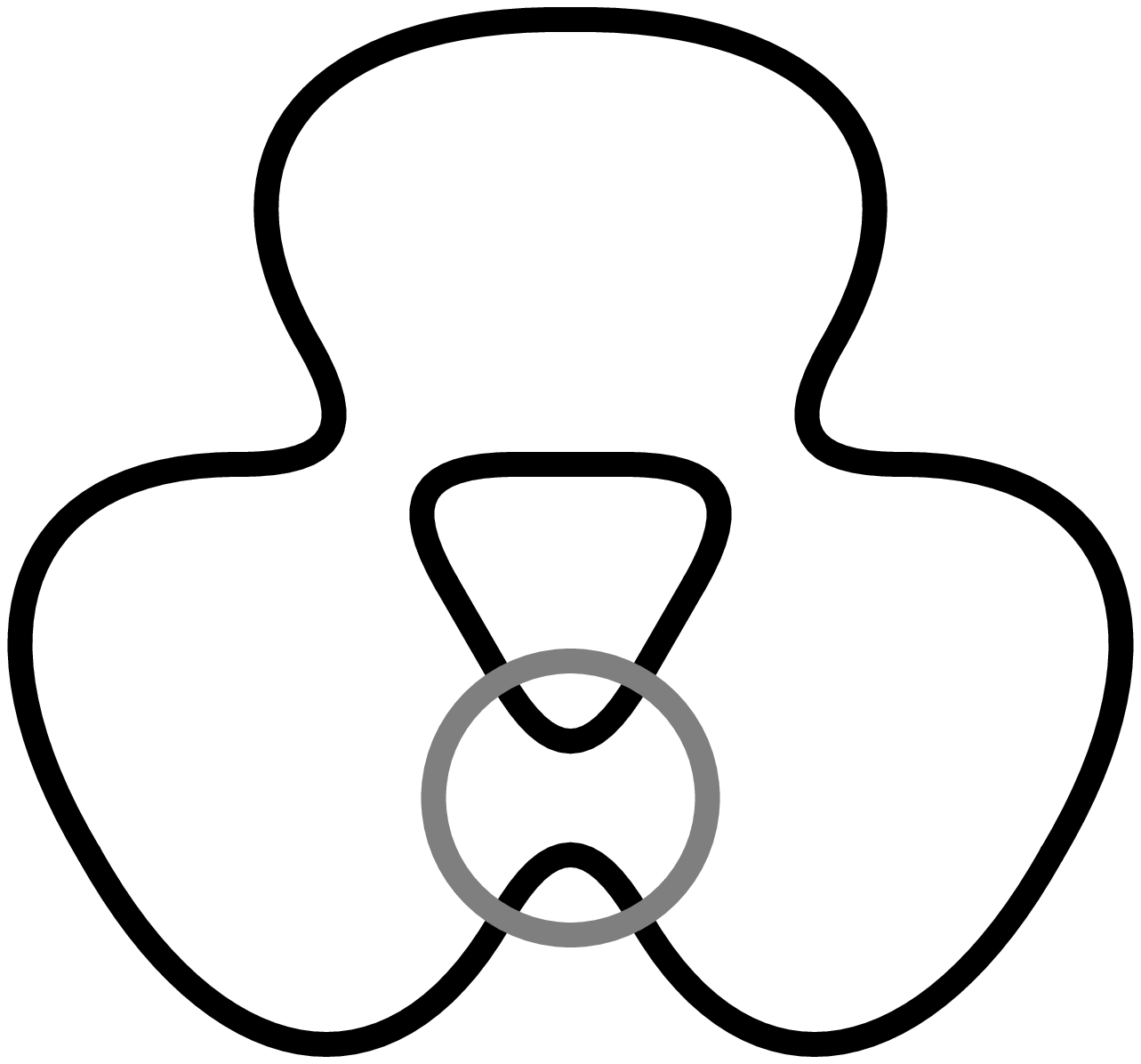}\\[-5pt] \textrm{\tiny{(0,0,0)}}\end{array}}\ar[uurr]^-{m}\ar[rr]^-{m}\ar[ddrr]_-{m}&& {\begin{array}{c}\includegraphics[scale=0.09375]{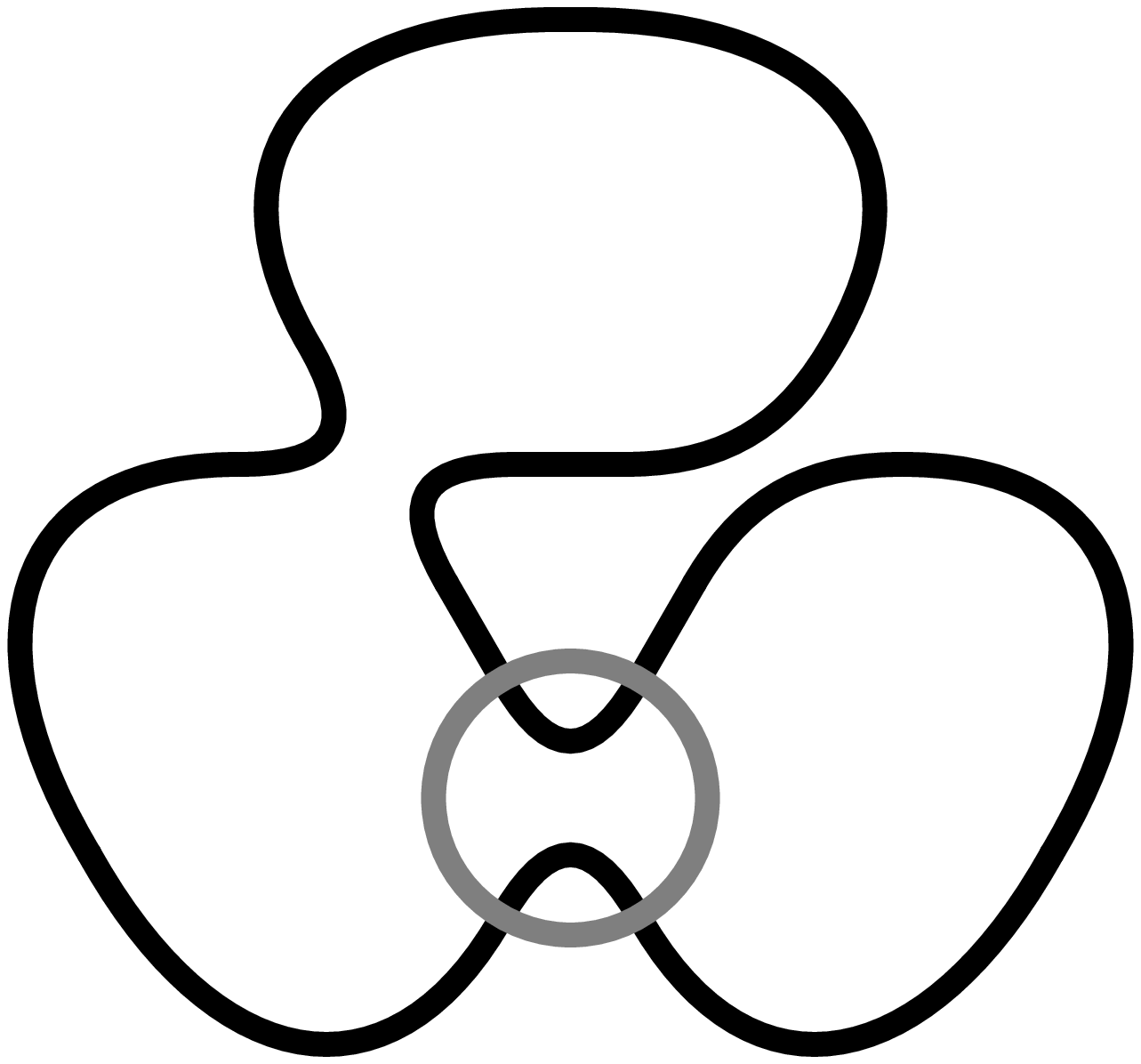}\\[-5pt]\textrm{\tiny{(0,1,0)}}\end{array}}\ar@{..}[d]\ar[uurr]^(0.65){\Delta}\ar[ddrr]_(0.65){-\Delta} & &{\begin{array}{c}\includegraphics[scale=0.09375]{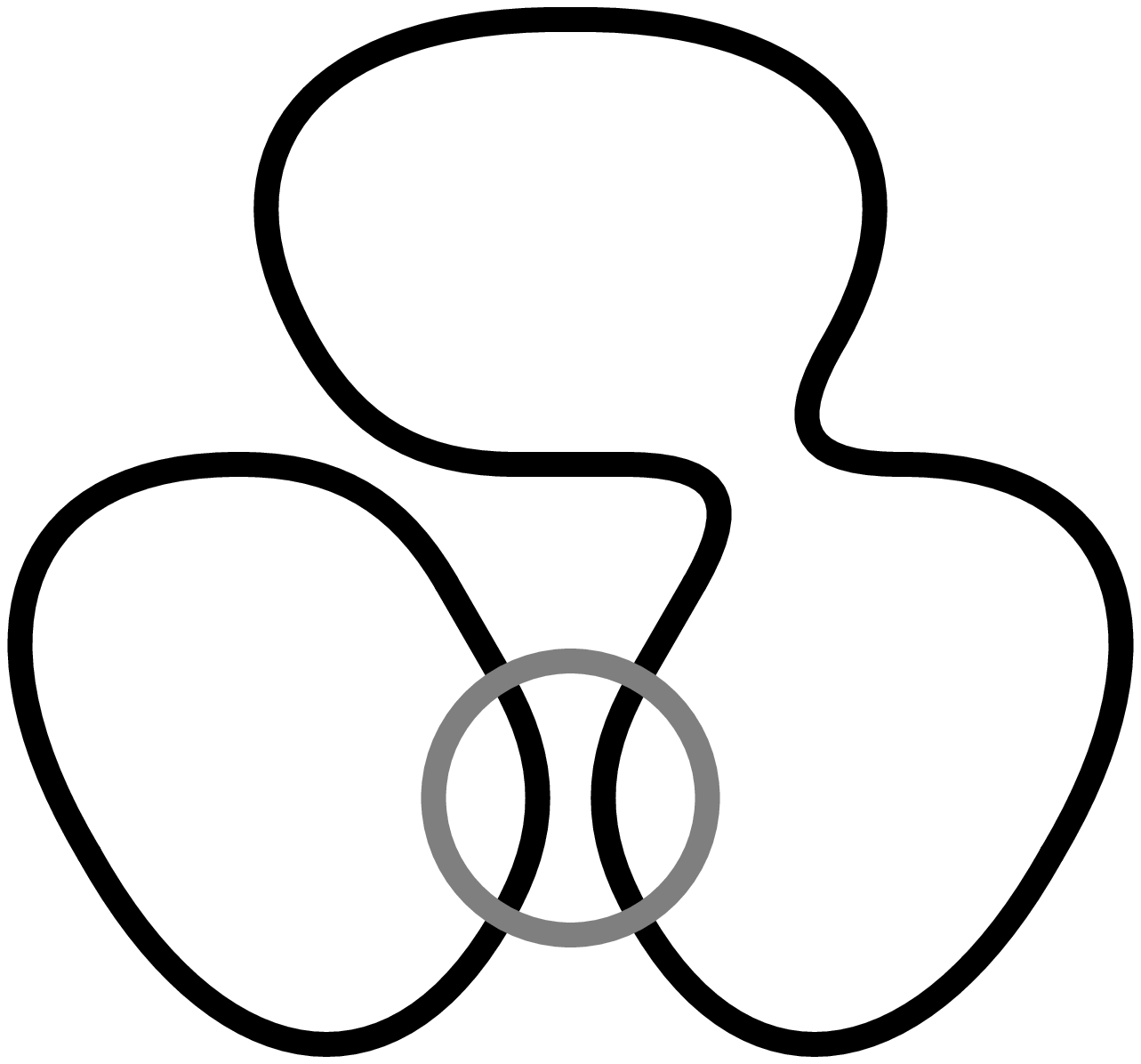}\\[-5pt]\textrm{\tiny{(1,0,1)}}\end{array}}\ar@{..}[d]\ar[rr]^-{-\Delta} & &{\begin{array}{c}\includegraphics[scale=0.09375]{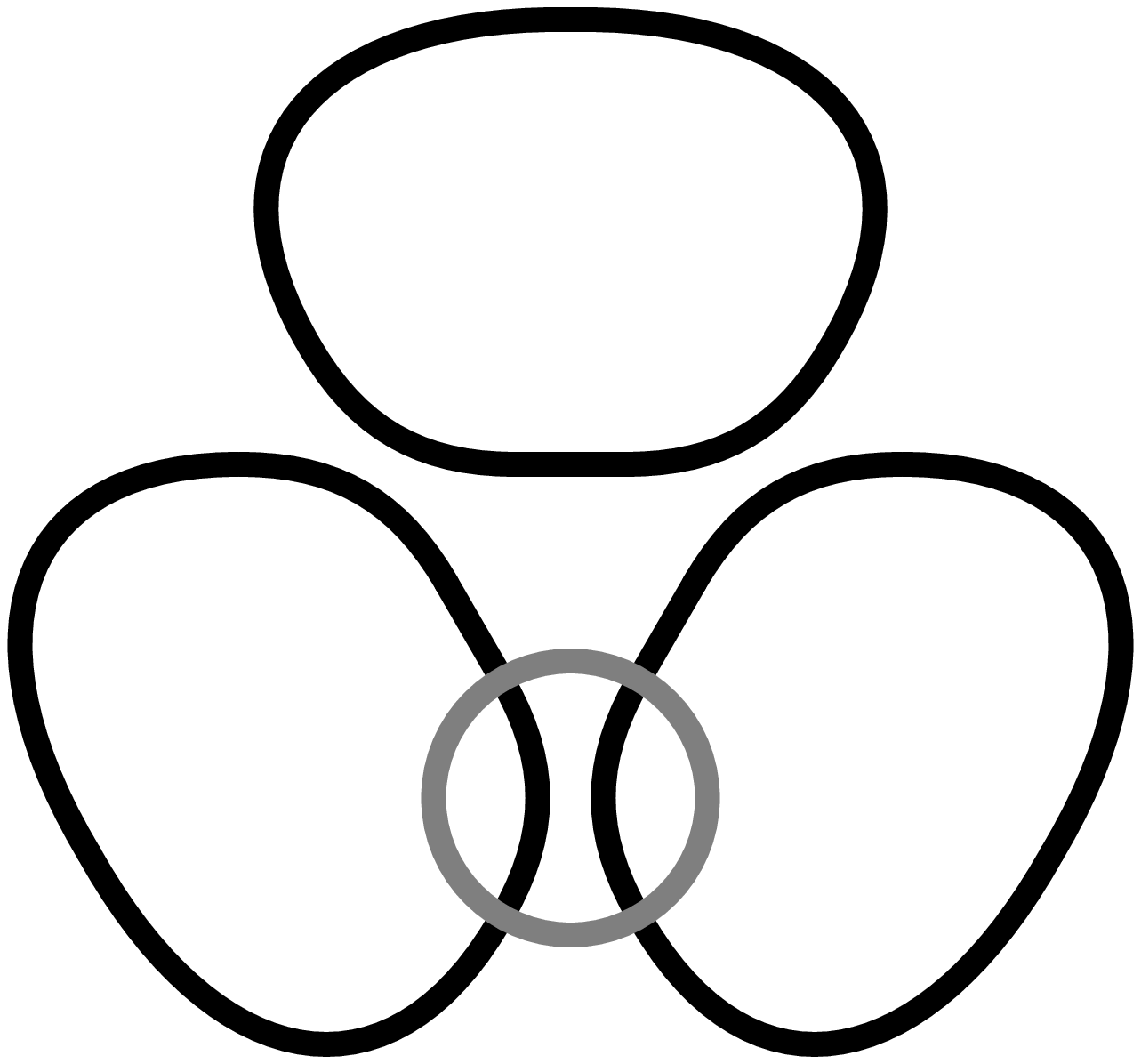}\\[-5pt]\textrm{\tiny{(1,1,1)}}\end{array}}\\
  && {\oplus}\ar@{..}[d] & \ar[ur]^(0.3){\Delta} & {\oplus}\ar@{..}[d] && \\
 & &{\begin{array}{c}\includegraphics[scale=0.09375]{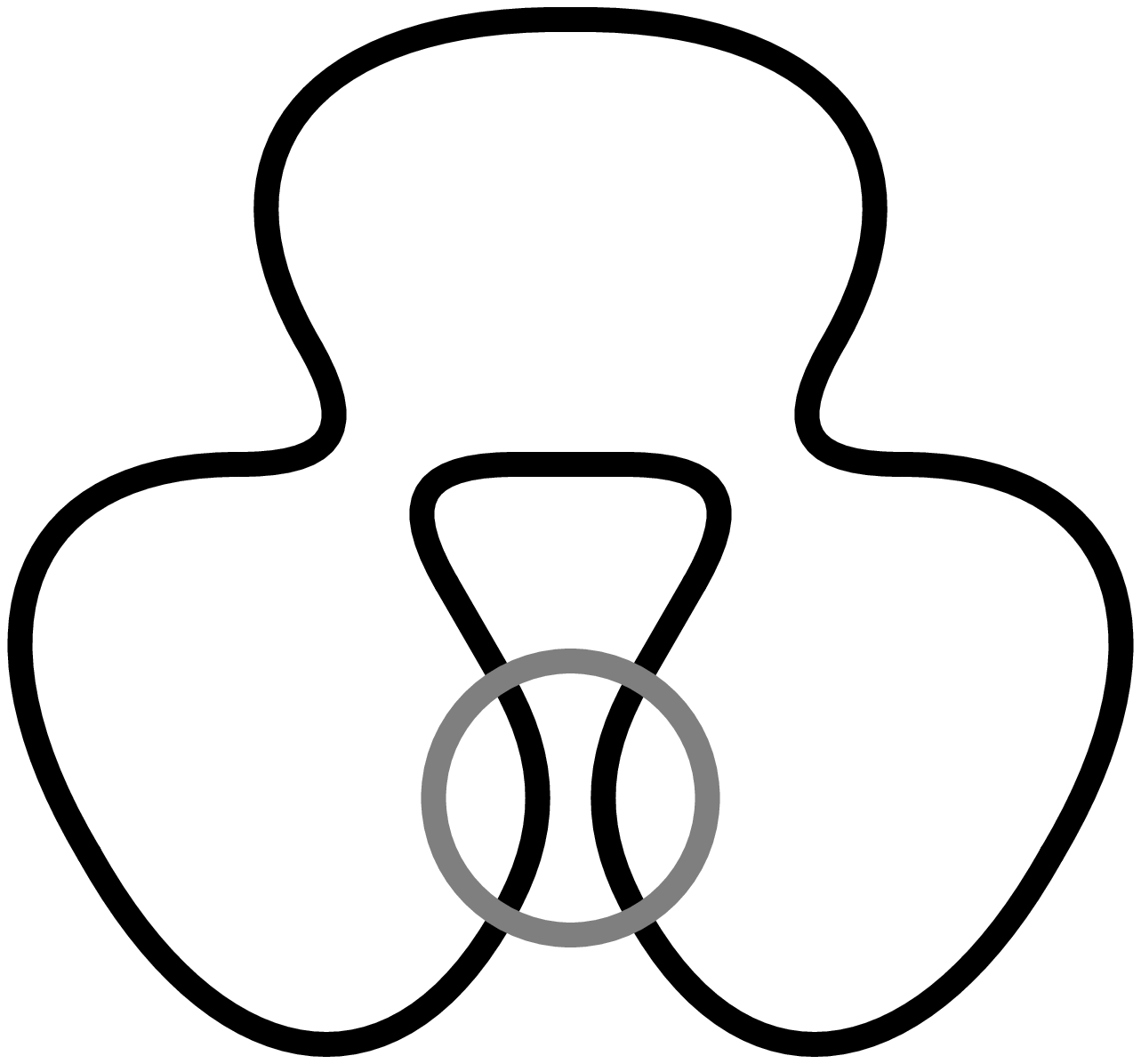} \\[-5pt] \textrm{\tiny{(0,0,1)}} \end{array}}\ar@{-}[ur]\ar[rr]_-{\Delta} && {\begin{array}{c}\includegraphics[scale=0.09375]{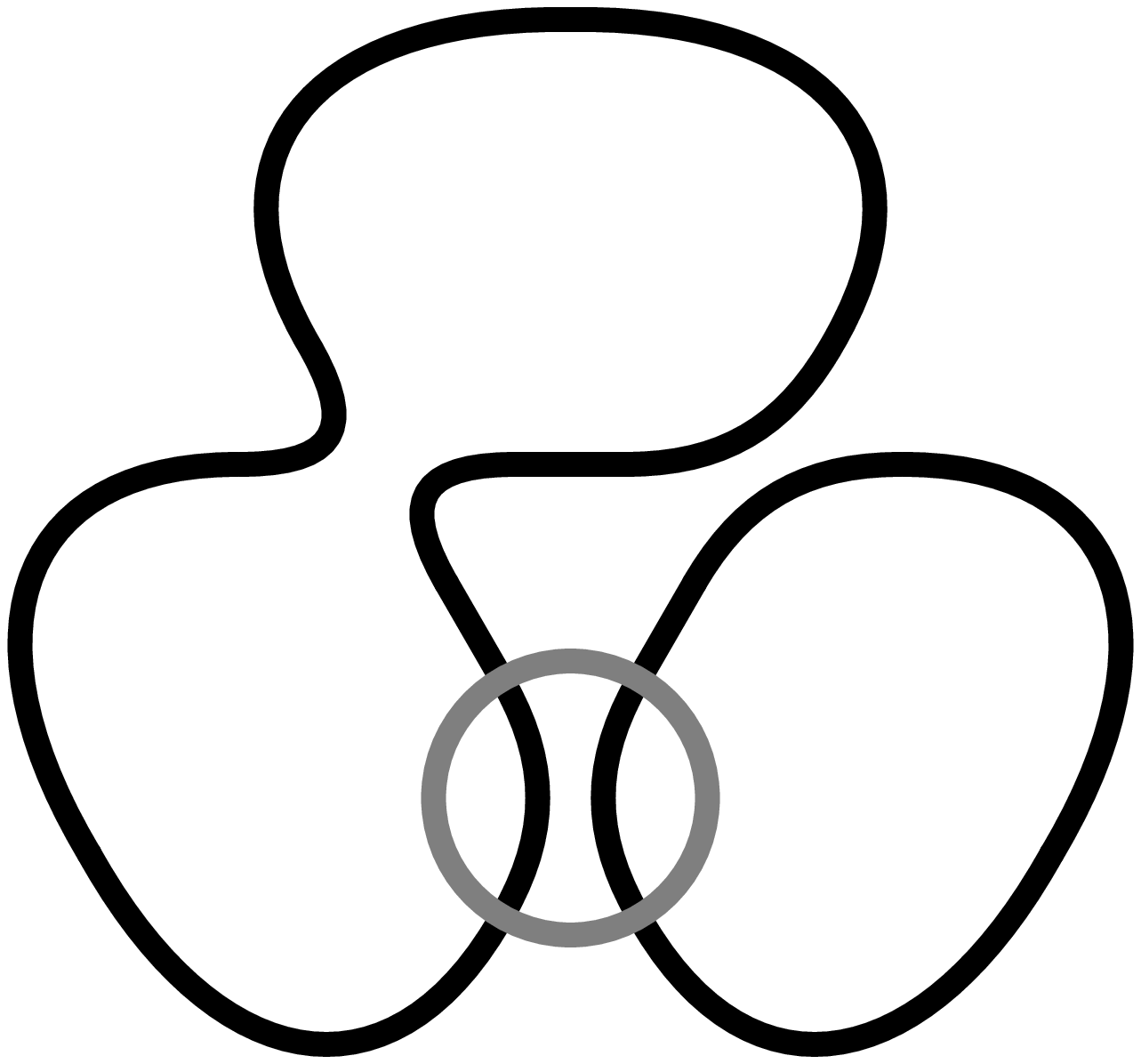} \\[-5pt] \textrm{\tiny{(0,1,1)}} \end{array}}\ar[uurr]_-{\Delta} & &
}
\end{array}$}
We briefly review Khovanov homology to solidify notation, and refer the reader to Khovanov's original paper \cite{Khovanov2000}, as well as \cite{Bar-Natan2002,Bar-Natan2005,Rasmussen2005,Turner2006}.

The Khovanov complex of a knot $K$ is generated by first considering an $n$-crossing diagram for $K$ together with $2^n$ states, each of which is a collection of disjoint simple closed curves in the plane. Each state $s$ is obtained from a choice of resolution $\zero$ (the 0-resolution) or $\one$ (the 1-resolution) for each crossing $\positive$. As a result, each state $s$ may be represented by an $n$-tuple with entries in $\{0,1\}$ so that the states may be arranged at the vertices of the $n$-cube $[0,1]^n$ (the cube of resolutions for $K$). Let $|s|$ be the sum of the entries of the $n$-tuple associated to $s$ (the height of $s$).  

Let $V$ be a free, graded $\bZ$-module generated by $\langle v_-,v_+\rangle$, where $\deg(v_\pm)=\pm1$. To each state we associate 
$V^{\ten \ell_s}$ where $\ell_s>0$ is the number of closed curves in the given state. The associated grading is referred to as the Jones grading, denoted by $q$. Set 
\[\sC^u(K) = \bigoplus_{u=|s|} V^{\ten \ell_s}\{|s|\}\]
where $\{\cdot\}$ shifts the Jones grading via $(W\{j\})_q = W_{q-j}$. The chain groups of the Khovanov complex are given by 
\[CKh^u_q(K) = \left(\sC(K)[-n_-]\{n_+-2n_-\}\right)^u_q = \sC^{u+n_-}_{q-n_++2n_-}(K)\]
where $[\cdot]$ shifts the homological grading $u$ as shown. For a given orientation of $K$, $n_+=n_+(K)$ is the number of positive crossings $\otherrightcross$ in $K$ and $n_-=n_-(K)$ is the number of negative crossings $\otherleftcross$ in $K$. The writhe of $K$ is given by $w(K)=n_+(K)-n_-(K)$.

The differentials $\partial^u : CKh^u(K)\to CKh^{u+1}$ come from the collection of edges in the cube of resolutions moving from height $u$ to height $u+1$. Each of these edges corresponds to exactly one of two operations 
\[
\begin{array}{cc}
\begin{aligned}
m\  :\ \ \ V \ten V&\longto V  \\[-4pt]
\vmm&\longmap 0 \\[-4pt]
\vp\ten v_\pm&\longmap v_\pm
\end{aligned} 
\quad
&
\quad
\begin{aligned}
\Del  \ :\ \ \ V&\longto V\ten V \\[-4pt]
\vm&\longmap\vmm\\[-4pt]
\vp&\longmap\vmp+\vpm
\end{aligned}
\end{array}
\]
of a Frobenius algebra defined over $V$, since each edge can be identified with exactly one change of the form $\zero\to\one$ or $\one\to\zero$. Fixing a convention so that the faces of the cube anti-commute, $\partial^u$ is the sum of all the maps at the prescribed height. 

The Khovanov homology $Kh(K)$, defined as the homology of the complex $(CKh^u(K),\partial^u)$, is an invariant of the knot $K$. Defining the Poincar\'e polynomial \[\chi_{(u,q)}(K)= \displaystyle{\sum}_i \displaystyle{\sum}_j u^i q^j \dim(Kh_j^i(K)\ten\bQ) \]
we have that $\chi_{(-1,q)}(K)=\hat{J}(K)$, the unormalized Jones polynomial, with $\hat{J}(\unknot) = q^{-1}+q$. 

Given a knot $K(\otherrightcross)$ with a distinguished positive crossing, there is a short exact sequence 
\[0\longto \sC\left(K(\zero)\right)[1]\{1\}\longto \sC\left(K(\negative)\right) \longto \sC\left(K(\one)\right)\longto 0.\]
Since $K(\one)$ inherits the orientation of $K(\otherrightcross)$, we set $c=n_-\left(K(\zero)\right)-n_-\left(K(\otherrightcross)\right)$ for some choice of orientation on $K(\zero)$ to obtain
\[0\longto CKh^{u-c-1}_{q-3c-2}\left(K(\zero)\right)\longto CKh^u_q\left(K(\negative)\right) \longto CKh^u_{q-1}\left(K(\one)\right)\longto 0.\]
This short exact sequence gives rise to a long exact sequence
\[\cdots\longto Kh^{u-c-1}_{q-3c-2}\left(K(\zero)\right)\longto Kh^u_q\left(K(\negative)\right) \longto Kh^u_{q-1}\left(K(\one)\right)\stackrel{\delta_*}{\longto} Kh^{u-c}_{q-3c-2}\left(K(\zero)\right)\longto\cdots.\]

\parpic[r]{$\begin{array}{c}\includegraphics[scale=0.1875]{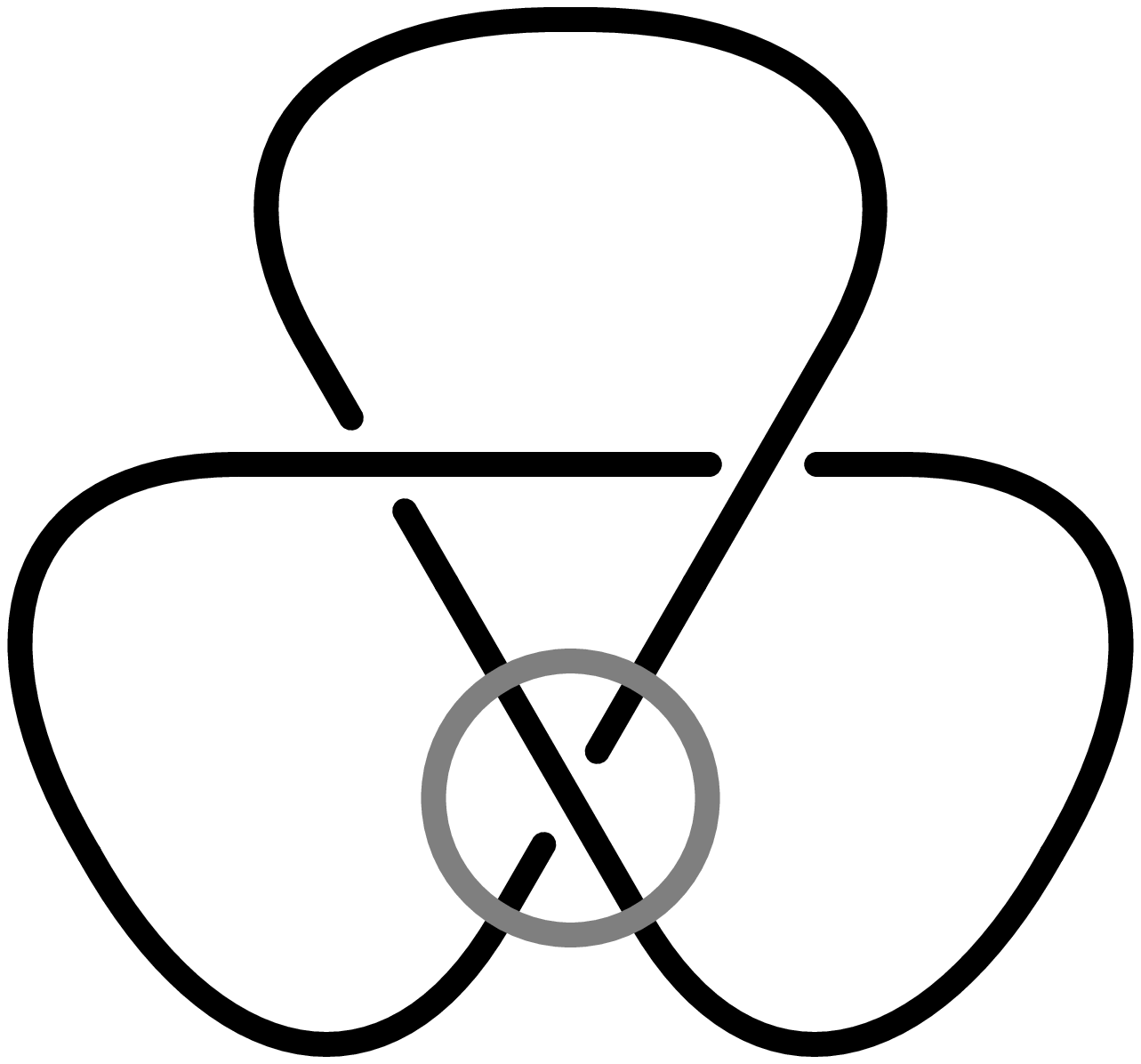}\end{array}$}
Here, $\delta_*$ is the map induced on homology from (the component of) the differential $\delta: CKh^u_{q-1}\left(K(\one)\right)\to CKh^{u-c-1}_{q-3c-2}\left(K(\zero)\right)$ in $CKh^u_{q-1}\left(K(\negative)\right)$. For example, in the complex for the right-hand trefoil (given on the previous page) we have circled the distinguished positive crossing (it is shown on the right). One can see that the subcomplex is given by states of the form $(\star,\star,1)$, and $\delta$ is given by maps that take states of the form $(\star,\star,0)$ to $(\star,\star,1)$.

Similarly, for a knot $K(\otherleftcross)$ with a distinguished negative crossing there is a long exact sequence
\[\cdots\longto Kh^{u}_{q+1}\left(K(\one)\right)\longto Kh^u_q\left(K(\positive)\right) \longto Kh^{u-c}_{q-3c-1}\left(K(\zero)\right)\stackrel{\delta_*}{\longto} Kh^{u+1}_{q+1}\left(K(\one)\right)\longto\cdots.\] 

\section{A Particular Example}\label{Particular Example} 

\parpic[r]{$\begin{array}{c}\includegraphics[scale=0.1875]{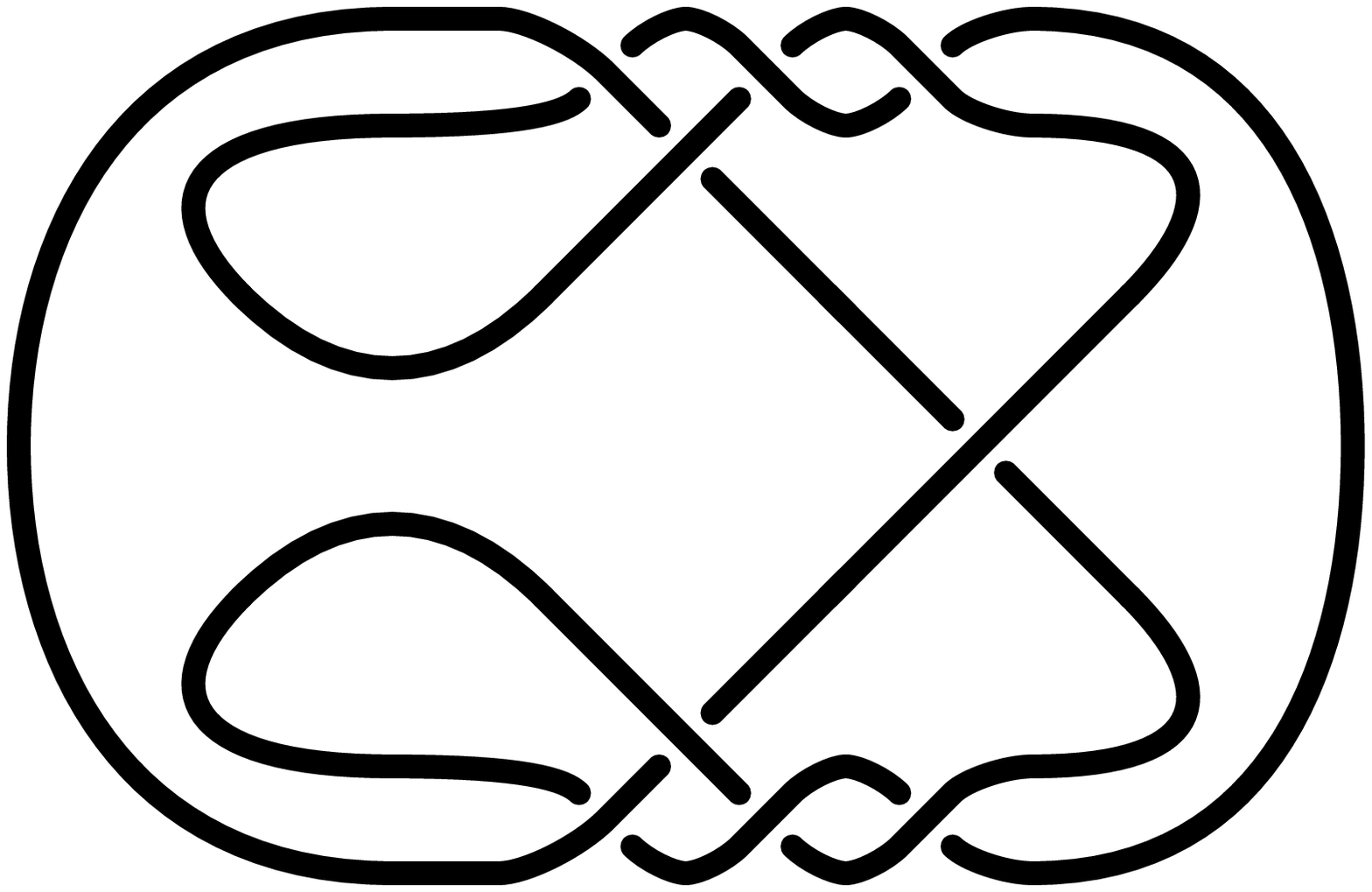}\qquad\includegraphics[scale=0.1875]{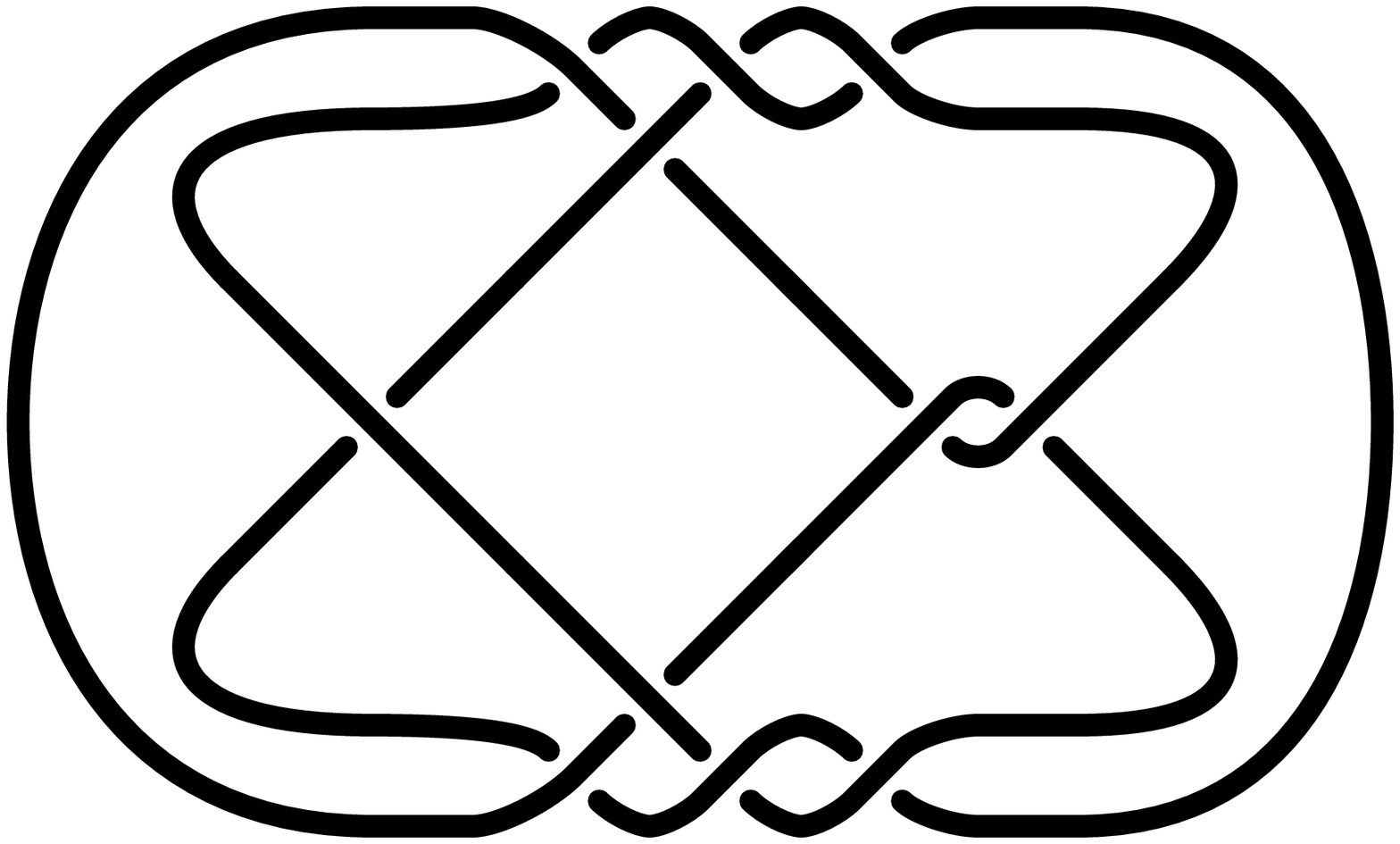}\end{array}$}
The knots $8_8$ and $10_{129}$ are shown -- up to mirrors -- on the right. These knots are known to have the same Jones polynomial, and in fact, they also have isomorphic Khovanov homology. The aim of this section is to prove this second fact making use of the long exact sequences from section \ref{Notation}. The calculation is made possible by this choice of (non-standard) projection for each of these knots. Indeed, from these projections we can see that the two knots are related by two local changes to the diagrams. We denote $10_{129} = K(\positive\negative)$ to distinguish the two (circled) crossings:
\[K(\positive\negative)=\raisebox{-0.67in}{\includegraphics[scale=0.25]{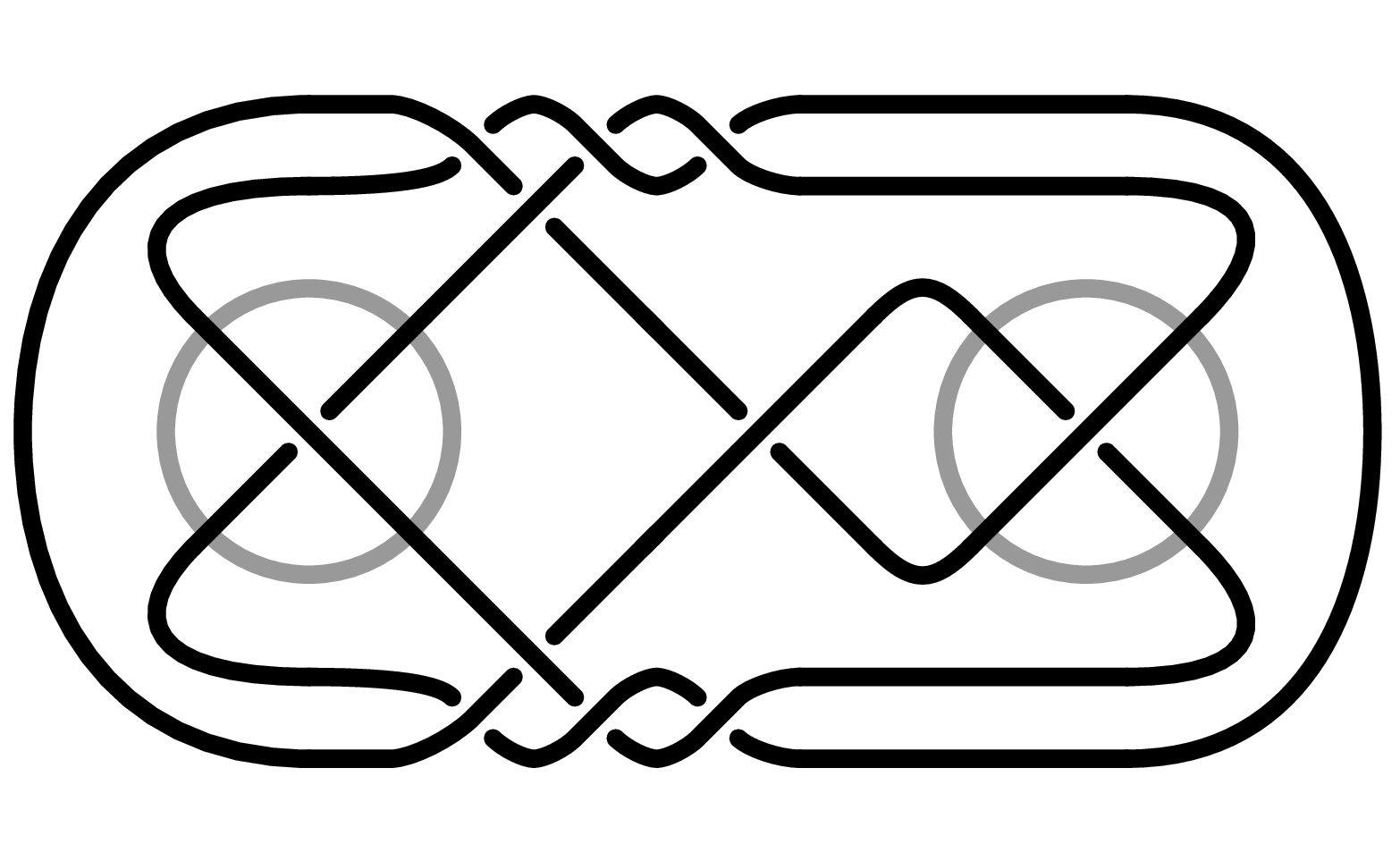}}\]    
With this notation, we have that $8_8=K(\zero\zero)$, and it can be checked that $K(\one\negative)$ and $K(\zero\one)$ (as well as $K(\negative\one)$ and $K(\one\zero)$) are diagrams for the two component trivial link. Let $L=Kh(\unknot\unknot)=V\ten V \cong (\bZ)_{-2}\pls(\bZ\pls\bZ)_0\pls(\bZ)_2$ and note that $L^u=0$ in all homological degrees $u\ne 0$. 

We begin by resolving with respect to the left-most crossing. Since this crossing is negative (choosing some orientation), we have the exact sequence
\[L_{q+1}^u\longto Kh^u_q K(\positive\negative) \longto Kh^{u+1}_{q+2} K(\zero\negative) \longto L^{u+1}_{q+1}\]
where one can check that $c=n_-(K(\zero\negative))-n_-(K(\positive\negative)) = -1$. In particular, for all homological gradings $u>0$, $L^u=0$ and 
\[  Kh^u_q K(\positive\negative) \cong Kh^{u+1}_{q+2} K(\zero\negative).\] 

Now we can resolve the remaining (distinguished) crossing in $K(\zero\negative)$. This crossing is positive, so we obtain a second exact sequence
\[L_{q-1}^{u-1}\longto Kh^{u-1}_{q-2}K(\zero\zero) \longto Kh^u_qK(\zero\negative) \longto L^u_{q-1}\]
where this time $c=n_-(K(\zero\zero))-n_-(K(\zero\negative)) = 0$. We can rewrite this as
\[L_{q+1}^{u}\longto Kh^{u}_{q}K(\zero\zero) \longto Kh^{u+1}_{q+2}K(\zero\negative) \longto L^{u+1}_{q+1}\]
so that for all homological gradings $u>0$, $L^u=0$ and 
\[Kh^{u}_{q}K(\zero\zero) \cong Kh^{u+1}_{q+2}K(\zero\negative).\]

In particular we have constructed a diagram of exact sequences
\[\xymatrix@R=7pt@C=30pt{
{L_{q+1}^{u}}\ar[r] & {Kh^{u}_{q}K(\zero\zero)}\ar@/^1pc/[dr] & & \\
{L_{q+1}^{u}}\ar[r] & {Kh^u_q K(\positive\negative)}\ar[r] & 	{Kh^{u+1}_{q+2}K(\zero\negative)}\ar[r]\ar@/_1pc/[dr] & {L^{u+1}_{q+1}} \\
& & & {L^{u+1}_{q+1}}
}\]
giving rise to an isomorphism
\[Kh^u_q K(\positive\negative) \cong Kh^{u}_{q}K(\zero\zero)\]
whenever $u>0$. Moreover, when $u=0$ we have that 
\[\xymatrix@R=7pt@C=30pt{
{(\bZ)_{-2}\pls(\bZ\pls\bZ)_0\pls(\bZ)_2}\ar[r] & {Kh^{0}_{q}K(\zero\zero)}\ar@/^1pc/[dr] & & \\
{(\bZ)_{-2}\pls(\bZ\pls\bZ)_0\pls(\bZ)_2}\ar[r] & {Kh^0_q K(\positive\negative)}\ar[r] & 	{Kh^1_{q+2}K(\zero\negative)}\ar[r]\ar@/_1pc/[dr] & {0} \\
& & & {0}
}\]
so that 
\[Kh^0_q K(\positive\negative) \cong Kh^0_{q}K(\zero\zero)\]
for $q\ne-3,-1,1$.

On the other hand, we may construct a second such diagram by considering the long exact sequence obtained by smoothing the (distinguished) right-most crossing of $K(\positive\negative)$ first. This gives the exact sequence
\[L_{q-1}^{u-1}\longto Kh^{u-1}_{q-2}K(\positive\zero) \longto Kh^u_qK(\positive\negative) \longto L^u_{q-1}\]
since the crossing we are resolving is positive, and $c=n_-(K(\positive\negative))-n_-(K(\positive\zero))=0$. This time we can conclude that 
\[Kh^{u-1}_{q-2}K(\positive\zero) \cong Kh^u_qK(\positive\negative)\]
whenever $u<0$. Resolving with respect to the remaining (distinguished) crossing of $K(\positive\zero)$ gives 
\[L_{q-1}^{u-1}\to Kh^{u-1}_{q-2}K(\positive\zero) \to Kh^u_qK(\zero\zero )\to L^u_{q-1}\]
(where we have shifted the gradings accordingly as before) since the crossing is negative, and $c=n_-(K(\zero\zero ))-n_-(K(\positive\zero))=-1$. Thus 
\[Kh^{u-1}_{q-2}K(\positive\zero) \cong Kh^u_qK(\zero\zero )\]
whenever $u<0$. This time the diagram of exact sequences is
\[\xymatrix@R=7pt@C=30pt{
	{L_{q-1}^{u-1}}\ar@/^1pc/[dr] &  & & \\
{L_{q-1}^{u-1}}\ar[r] & {Kh^{u-1}_{q-2}K(\positive\zero)}\ar[r]\ar@/_1pc/[dr] & {Kh^u_q K(\positive\negative)}\ar[r] & {L^{u}_{q-1}} \\
& & {Kh^u_qK(\zero\zero )}\ar[r] & {L^{u}_{q-1}}
}\]
giving rise to an isomorphism 
\[Kh^u_q K(\positive\negative) \cong Kh^{u}_{q}K(\zero\zero)\]
whenever $u<0$. When $u=0$, we have that 
\[\xymatrix@R=7pt@C=30pt{
	{0}\ar@/^1pc/[dr] &  & & \\
{0}\ar[r] & {Kh^{-1}_{q-2}K(\positive\zero)}\ar[r]\ar@/_1pc/[dr] & {Kh^0_q K(\positive\negative)}\ar[r] & {(\bZ)_{-2}\pls(\bZ\pls\bZ)_0\pls(\bZ)_2} \\
& & {Kh^0_qK(\zero\zero )}\ar[r] & {(\bZ)_{-2}\pls(\bZ\pls\bZ)_0\pls(\bZ)_2}
}\]
so that
\[Kh^0_q K(\positive\negative) \cong Kh^0_{q}K(\zero\zero)\]
for $q\ne-1,1,3$.

Combining the information gained from both diagrams, we can conclude that
\[Kh^u_q K(\positive\negative) \cong Kh^{u}_{q}K(\zero\zero)\]
except when $u=0$ and $q=\pm1$. In fact, the above diagram tells us that the torsion parts for $u=0$ and $q=\pm1$ are isomorphic. Indeed, since $L$ is torsion free we have 
\[\xymatrix@R=7pt@C=30pt{
	{0}\ar@/^1pc/[dr] &  & & \\
{0}\ar[r] & {\left(Kh^{-1}_{q-2}K(\positive\zero)\right)_{\operatorname{tor}}}\ar[r]\ar@/_1pc/[dr] & {\left(Kh^0_q K(\positive\negative)\right)_{\operatorname{tor}}}\ar[r] & {0} \\
& & {\left(Kh^0_qK(\zero\zero )\right)_{\operatorname{tor}}}\ar[r] & {0}
}\]
hence 
\[\left(Kh^0_q K(\positive\negative)\right)_{\operatorname{tor}} \cong \left(Kh^0_q K(\zero\zero)\right)_{\operatorname{tor}}\]
for all $q$. On the other hand, the free part may be analyzed over $\bQ$. Note that since all groups are isomorphic away from $(u,q)=(0,\pm1)$ we have that
\begin{align*}
& \chi_{(-1,q)}(K(\positive\negative)) = \chi_{(-1,q)}(K(\zero\zero)) \\
\Longleftrightarrow & \displaystyle{\sum}_{j=\pm1} q^j\dim(Kh^0_jK(\positive\negative)\ten\bQ) = \displaystyle{\sum}_{j=\pm1} q^j\dim(Kh^0_jK(\zero\zero)\ten\bQ) \\
\Longleftrightarrow & \dim(Kh^0_{\pm1}K(\positive\negative)\ten\bQ) = \dim(Kh^0_{\pm1}K(\zero\zero)\ten\bQ) \\
\Longleftrightarrow & Kh^0_{\pm1}K(\positive\negative)\ten\bQ \cong Kh^0_{\pm1}K(\zero\zero)\ten\bQ.
\end{align*}
But since $\hat{J}(8_8)=\hat{J}(10_{129})$, we know that $\chi_{(-1,q)}(K(\positive\negative)) = \chi_{(-1,q)}(K(\zero\zero))$ and hence we get an isomorphism on the free part for $u=0$. Finally, we conclude that $Kh(8_8)\cong Kh(10_{129})$. Although this result was expected, the goal of sections \ref{Construction} and \ref{Proof} is to exploit this computation in a more general setting.

\section{Construction}\label{Construction}

\parpic[r]{$\begin{array}{c}\includegraphics[scale=0.25]{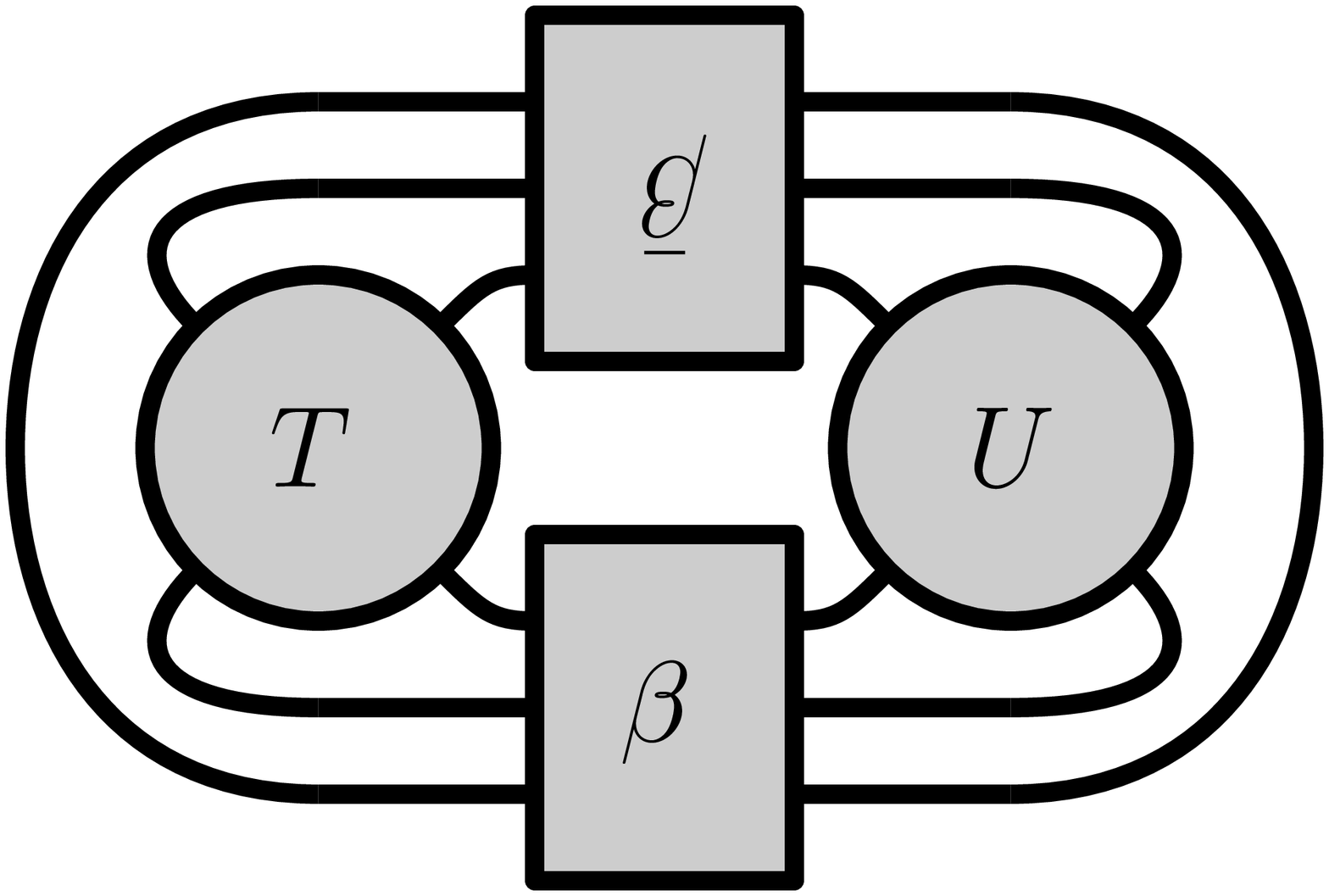}\end{array}$}
Consider the knot $K=K_\beta(T,U)$ (cf. \cite{Watson2005}) where $\be$ is an element of the three strand braid group $B_3$ with inverse $\bar{\be}$. $T$ and $U$ are tangles (or Conway tangles), that is $T=(B^3_T, \tau)$ and $U=(B^3_U, \mu)$ where $B^3_T$ (respectively $B^3_U$) is a 3-ball containing a collection of arcs $\tau$ (respectively $\mu$) that intersect the boundary of the 3-ball transversally in exactly 4 points \cite{Lickorish1981, Rolfsen1994}. 

There is a well defined $\bZ$-action (a half-twist action) on the set of isotopy classes (fixing endpoints) of tangles that comes from the 2-strand braid group. For a given tangle $T$, write 
 $T^\si=\sigright$ and $T^{\bar{\si}}=\SIGright$
 where $\langle\si\rangle=\bZ\cong B_2$ and $\si\bar{\si}=e$ (that is, $\si = \positive$ is the standard braid generator). Let $K^\si=K_\be\left(T^\si,U^{\bar{\si}}\right)$. The following is proved in \cite{Watson2005}.
\begin{lemma}\label{lem:jones}
If $w(K^\si)=w(K)$ then $\hat{J}(K^\si)=\hat{J}(K)$.
\end{lemma}

We now restrict our attention to braids of the form $\beta=\si_1^{-1}\si_2\si_1^{-2n}$ for simplicity, and assume that the tangles considered have no closed components (that is, $\tau$ and $\mu$ are each a pair of arcs). We require the tangles $T$ and $U$ of $K$ to compatible in the following sense: For a given orientation of $K$, the tangles $T$ and $U$ are compatible if $w(K_\be(T,U))=w(K_\be(T^\si,U^{\bar{\si}}))$ (this is equivalent to requiring that $n_\pm(K_\be(T,U))=n_\pm(K_\be(T^{\si},U^{\bar{\si}}))-1$). For example, if $U=T^\star$ where $T^\star$ is the mirror image of $T$, then U and T are compatible.

The sum $T+U$ of two tangles is defined by side-by-side concatenation. This generalizes the half-twist-action: $T^\si$ may be denoted $T + \Cpositive$ ($\si$ {\em adds} a twist). A tangle $T$ is called simple if $T+\Cone$ is isotopic (fixing endpoints) to $\Cone$. Note that $T^\si$ is simple iff $T$ is simple.

\begin{lemma}\label{lem:khovanov}
Suppose $(T,U)$ is a compatible, simple pair of tangles for $K_\beta(-,-)$. Then $K_\beta(T,U)$ and $K_\be\left(T^\si,U^{\bar{\si}}\right)$ have identical Khovanov homology.
\end{lemma} 

\begin{remark}\label{rem:left}
We could also consider an action that adds twists to tangles on the left; we make use of this in section \ref{Mutants}.
\end{remark}

\section{Proof}\label{Proof}

The proof of lemma \ref{lem:khovanov} is an application of the long exact sequences, and proceeds in much the same way as in the example of section \ref{Particular Example}. As in section \ref{Particular Example}, we distinguish two crossings of $K_\be\left(T^\si,U^{\bar{\si}}\right)$ (the two crossings added by the action of $\si$) and write $K_\beta(\positive\negative) = K_\be\left(T^\si,U^{\bar{\si}}\right)$ so that $K_\beta(\zero\zero)=K_\beta(T,U)$. Since the tangles $T$ and $U$ are simple by hypothesis, it is easy to check that $K_\beta(\one\negative)$ and $K_\beta(\zero\one)$ (as well as $K_\beta(\positive\one)$ and $K_\beta(\one\zero)$) are diagrams for the 2-component trivial link (note that the braids cancel). Again, we denote $L=Kh(\unknot\unknot)=V\ten V \cong (\bZ)_{-2}\pls(\bZ\pls\bZ)_0\pls(\bZ)_2$.

To exploit the calculation from section \ref{Particular Example}, we must verify that the required values of $c$ used in the long exact sequences follow from the compatibility hypothesis. First note that since we are considering knots, we may restrict attention to tangles $T$ and $U$ that have connectivity of the form $\Gcross$ and $\Gzero$ (tangles having connectivity $\Gone$ give rise to links). The connectivity of a tangle is simply the input and output data for each strand.

\parpic[r]{$\begin{array}{c}\includegraphics[scale=0.25]{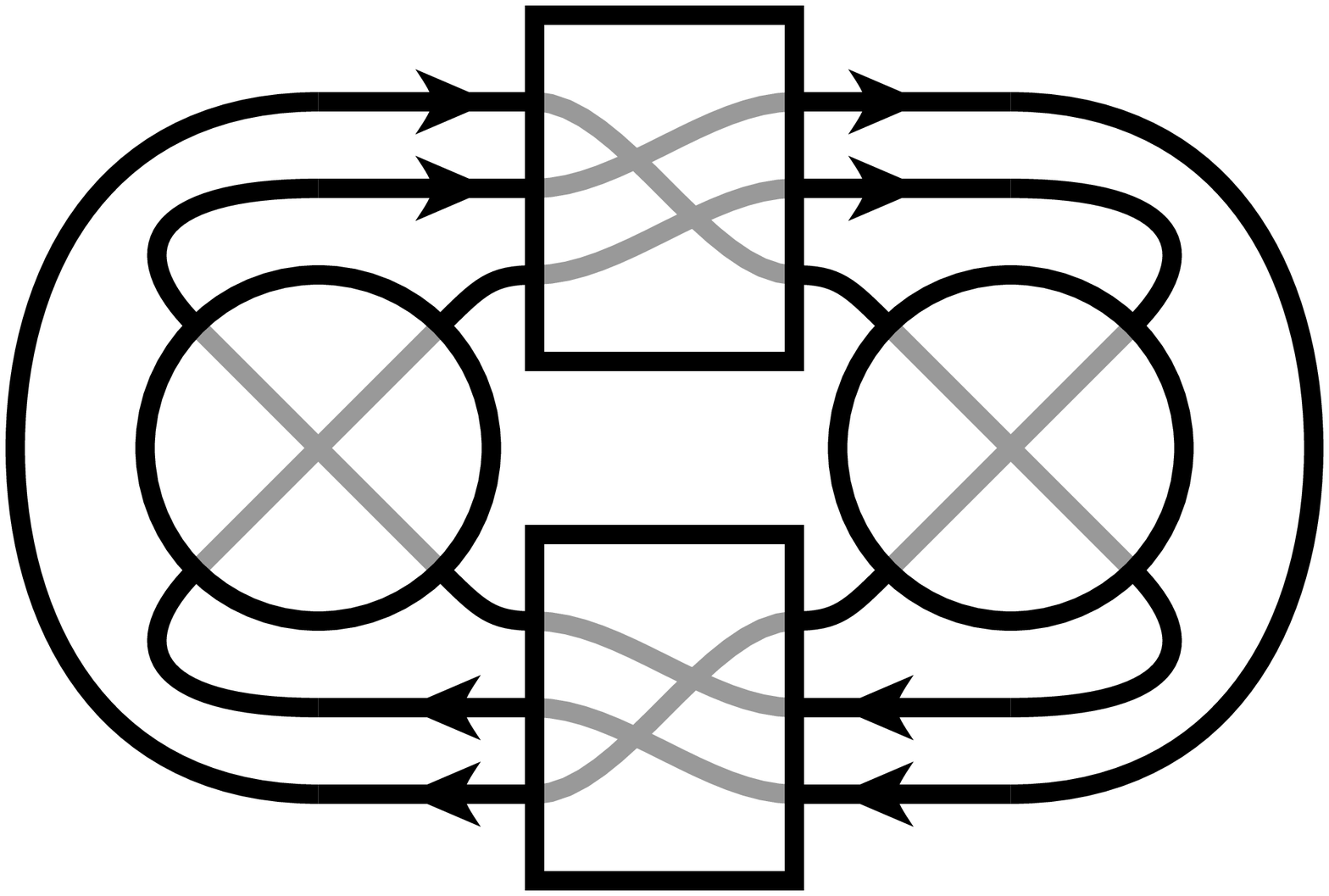}\end{array}$}
Suppose that both $T^\si$ and $U^{\bar{\si}}$ have connectivity of the form $\Gcross$. Then it is easy to check that $T$ and $U$ have connectivity of the form $\Gzero$. Moreover, since the permutation associated to $\beta=\si_1^{-1}\si_2\si_1^{-2n}$ is $(3\ 1\ 2)$ and the permutation associated to $\bar{\beta}=\si_1^{2n}\si_2^{-1}\si_1$ is $(3\ 2\ 1)$, we can fix the orientation for $K_\beta(\positive \negative)$ shown on the right. With this orientation in hand, we have that the distinguished crossing of $T^\si$ is negative, while the distinguished crossing of $U^{\bar{\si}}$ is positive: $K(\otherleftcross\otherrightcross)$. Notice that the resolution $K(\zero\negative)$ (that is, the tangle $T$) does not inherit this orientation. If we resolve with respect to the left-most (distinguished) crossing first, we have
\[\xymatrix@C=30pt{
{\includegraphics[scale=0.1875]{figures/skeleton.ps}}\ar[r]^{c=-1} &
{\includegraphics[scale=0.1875]{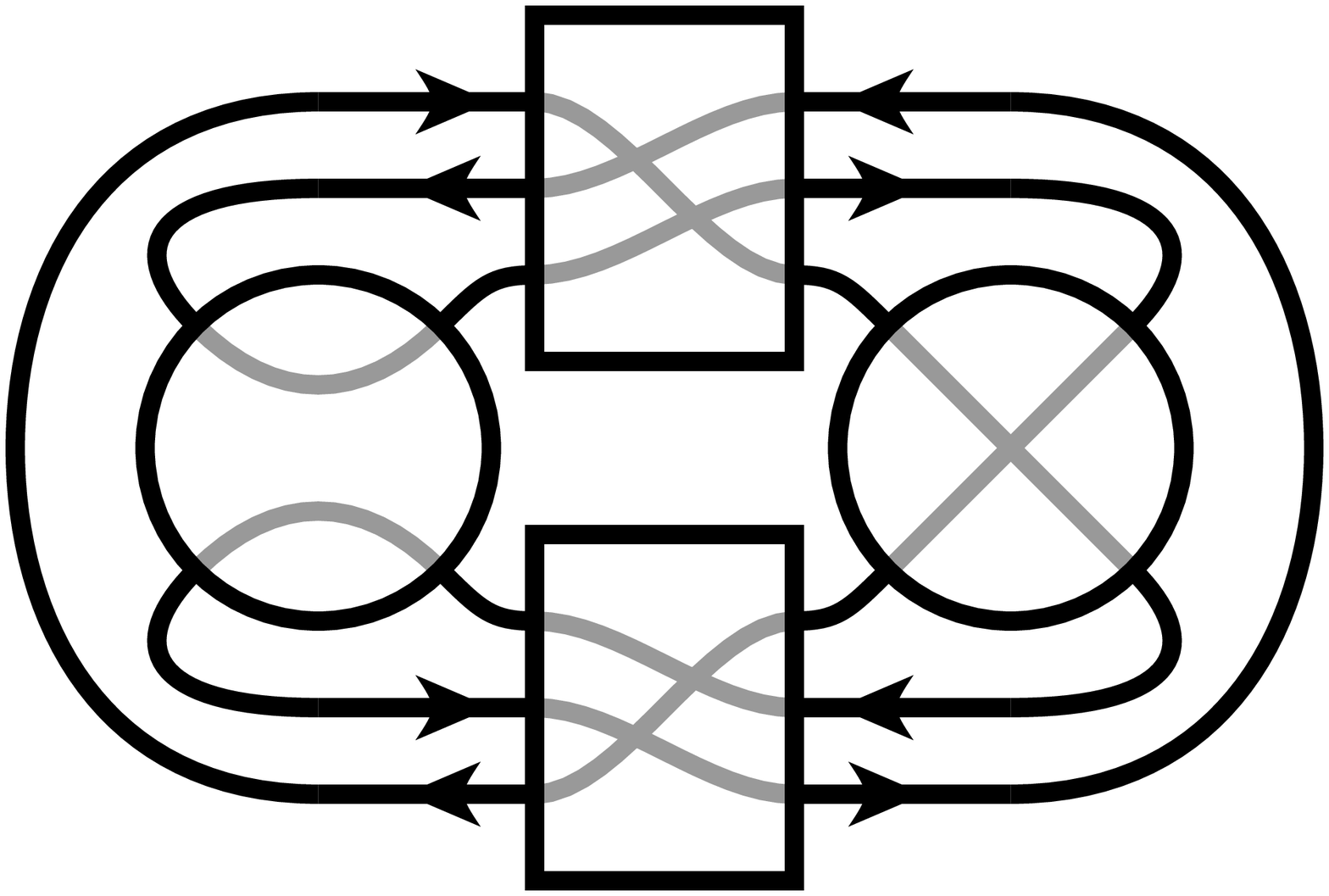}}\ar[r]^{c=0} &
{\includegraphics[scale=0.1875]{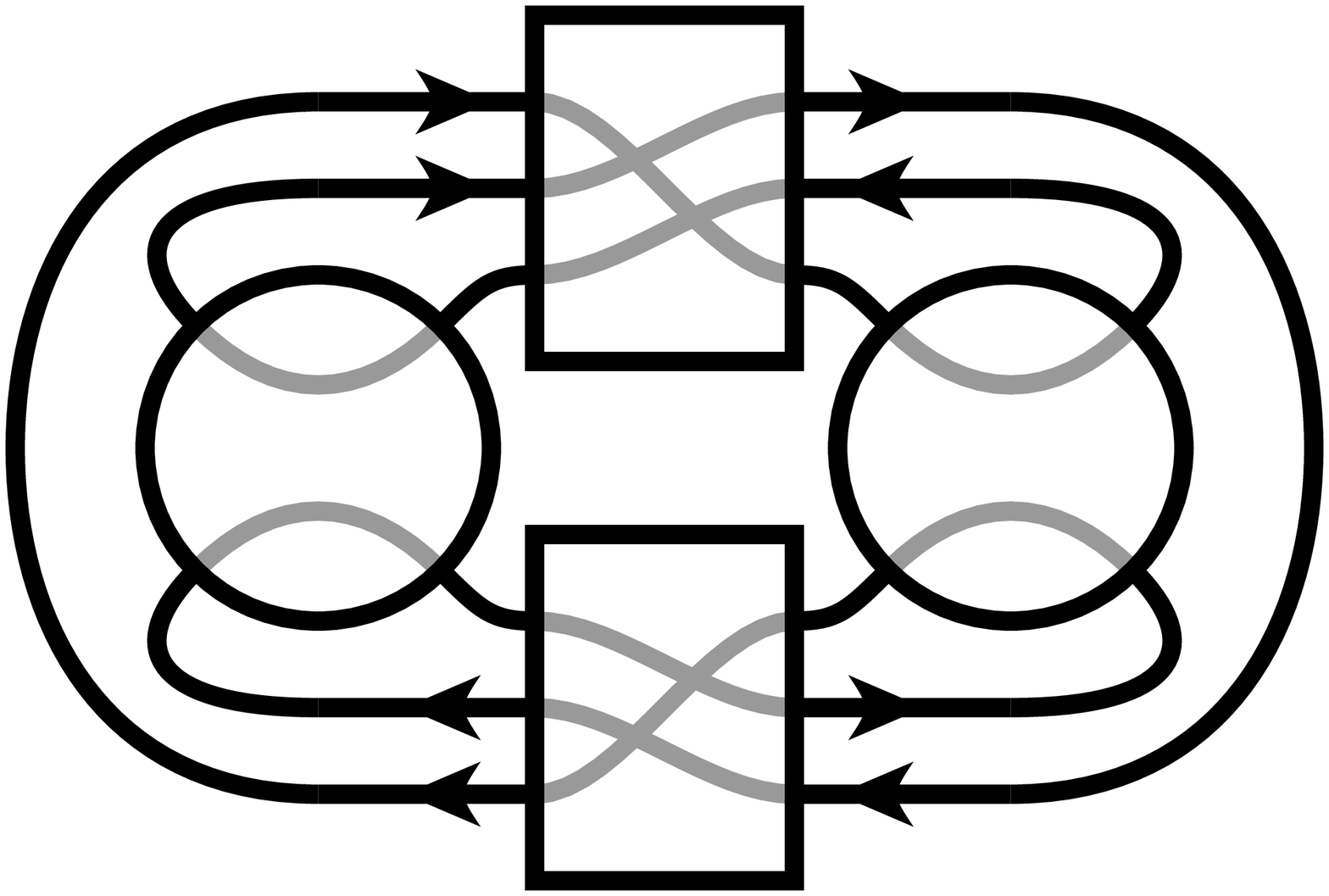}.}
}\]
\parpic[r]{$\begin{array}{c}\xymatrix@C25pt@R-25pt{ 
{\includegraphics[scale=0.1875]{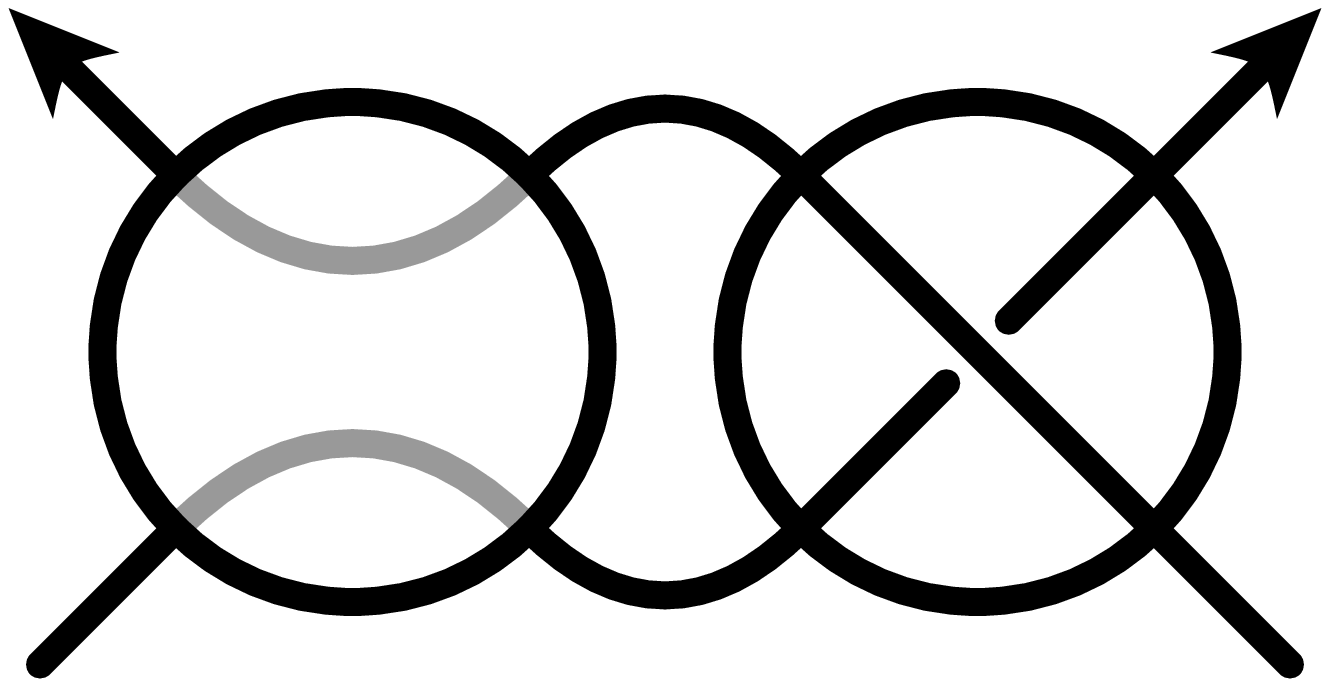}}\ar[r]^{c=-1} &
{\includegraphics[scale=0.1875]{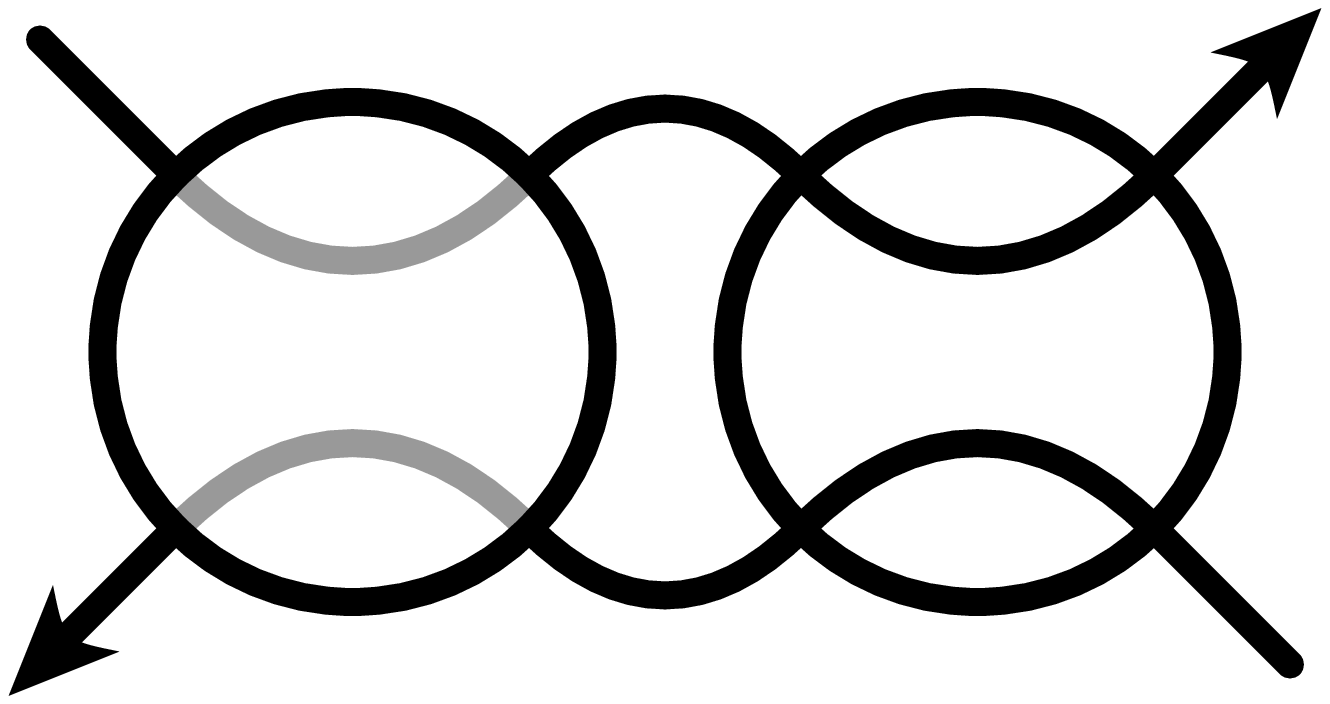}} \\
{\includegraphics[scale=0.1875]{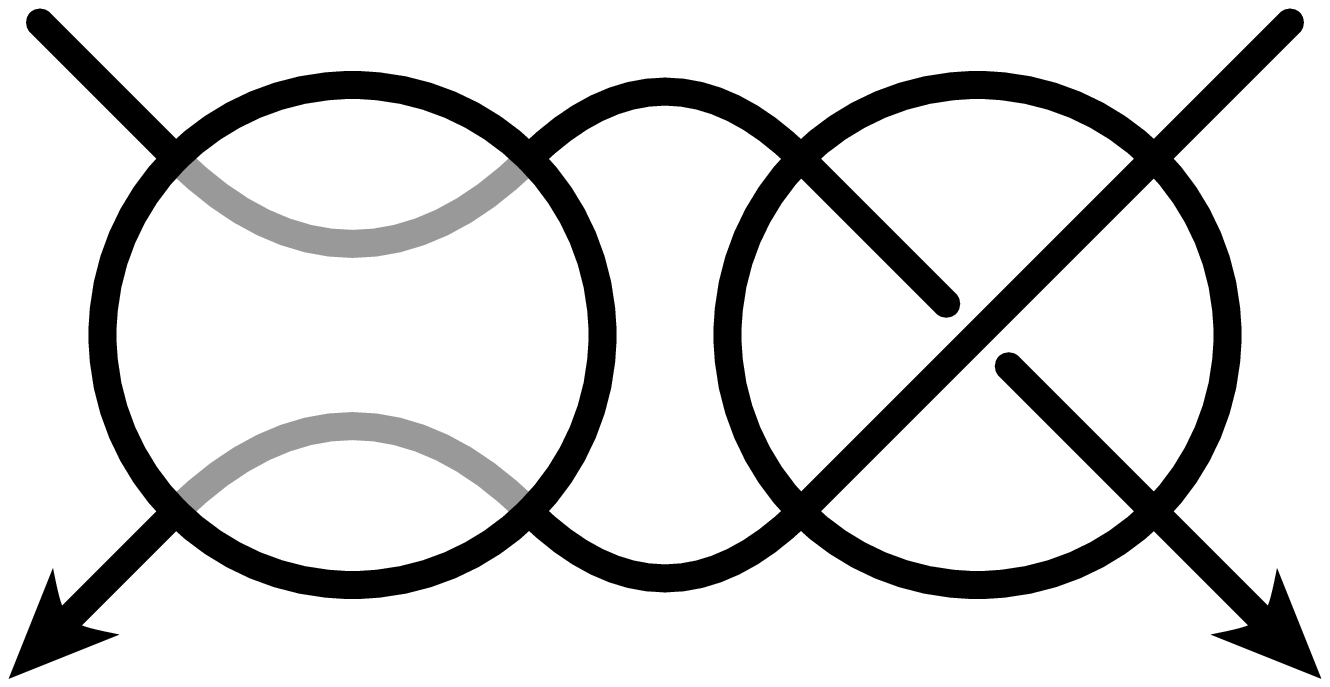}}\ar[r]^{c=0} &
{\includegraphics[scale=0.1875]{figures/resolution1.ps}}
}\end{array}$} To justify these values for $c$, we look at the resolutions at each stage. Of course, since the braids chosen are inverses of each other, the number of negative and positive crossings contributed by the braids remains constant. Therefore, to compute the values for $c$ we need only consider the tangles $T^\si$ and $U^{\bar{\si}}$ (shown on the right) upon resolution of the respected distinguished crossings. In the first step, the orientation on $U^{\bar{\si}}$ is preserved, while the new orientation for the resolution of $T^\si$ (that is, $T$) has precisely one less negative crossing (the crossing we resolved). This is because the new orientation reverses the orientation on both strands (one can check that this will always preserve the number of positive and negative crossings) so that $c=n_-(T)-n_-(T^\si)=n_-\left(K_\beta(\zero\negative)\right)-n_-\left(K_\beta(\positive\negative)\right)=-1$. Now, resolving the second crossing, we make similar observations. The orientation on the strands of the tangle $T$ are both reversed once more, so that the number of negative crossings contributed is left unchanged. On the other hand, he resolution taking $U^{\bar{\si}}$ to $U$ removes a positive crossing, and preserves the orientation on the tangle $U$. Therefore $c=n_-(U)-n_-(U^{\bar{\si}})=n_-\left(K_\beta(\zero\zero)\right)-n_-\left(K_\beta(\zero\negative)\right)=0$ as claimed.   

This, together with the previous observation that both $\one$-resolutions are unlinks, allows us to produce, exactly as in section \ref{Particular Example}, a diagram of long exact sequences:
\[\xymatrix@R=7pt@C=30pt{
{\cdots\ L_{q+1}^{u}}\ar[r] & {Kh^{u}_{q}K_\beta(\zero\zero)}\ar@/^1pc/[dr] & & &\\
{\cdots\ L_{q+1}^{u}}\ar[r] & {Kh^u_q K_\beta(\positive\negative)}\ar[r] & 	{Kh^{u+1}_{q+2}K_\beta(\zero\negative)}\ar[r]\ar@/_1pc/[dr] & {L^{u+1}_{q+1}\ \cdots} & {(1)}\\
& & & {L^{u+1}_{q+1}\ \cdots} &
}\]

On the other hand, resolving first with respect to the right-most (distinguished) crossing we have that
\[\xymatrix@C=30pt{
{\includegraphics[scale=0.1875]{figures/skeleton.ps}}\ar[r]^{c=0} &
{\includegraphics[scale=0.1875]{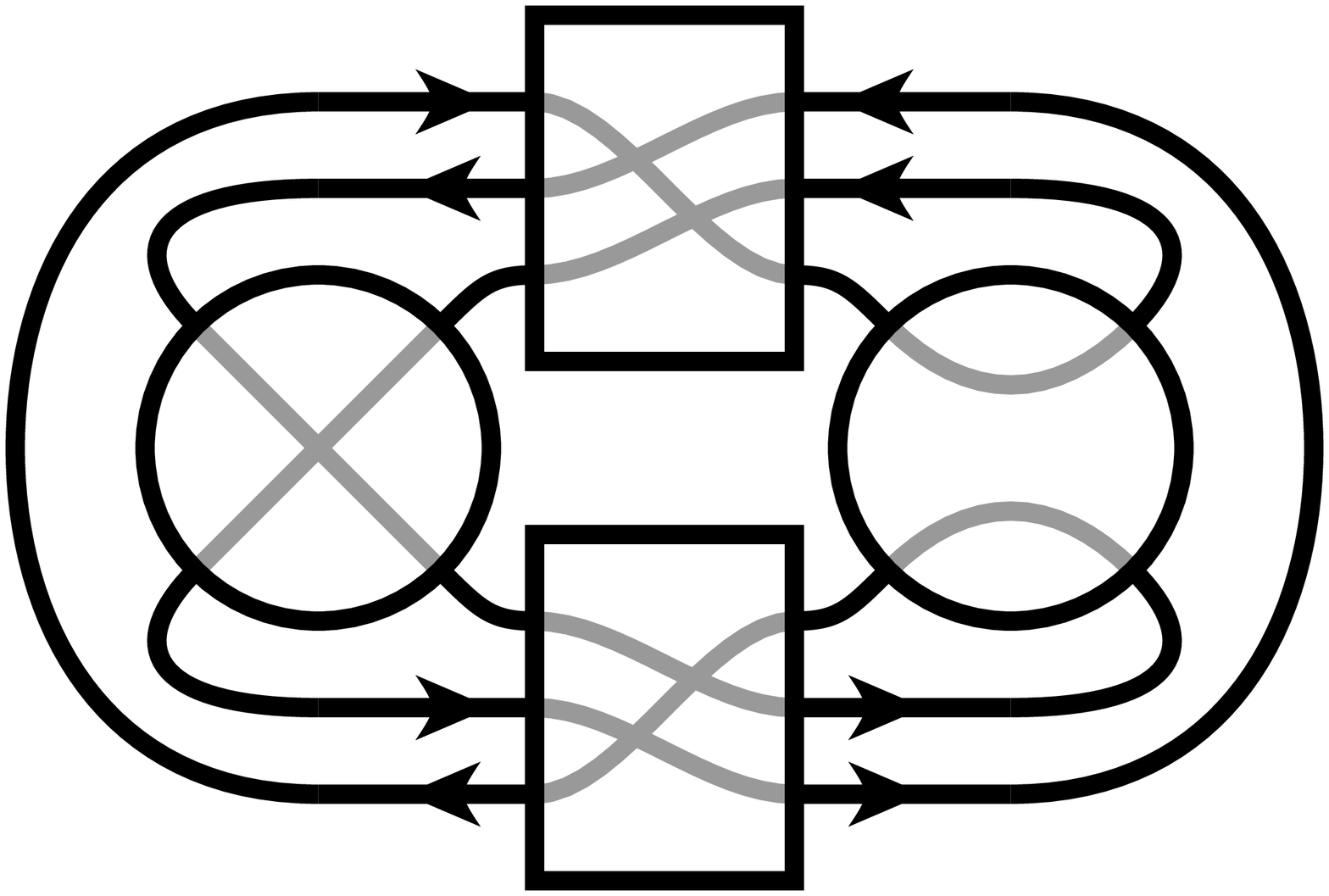}}\ar[r]^{c=-1} &
{\includegraphics[scale=0.1875]{figures/skeleton3.ps}}
}\]
\parpic[r]{$\begin{array}{c}\xymatrix@C25pt@R-25pt{ 
{\includegraphics[scale=0.1875]{figures/coherence2.ps}}\ar[r]^{c=0} &
{\includegraphics[scale=0.1875]{figures/resolution1.ps}} \\
{\includegraphics[scale=0.1875]{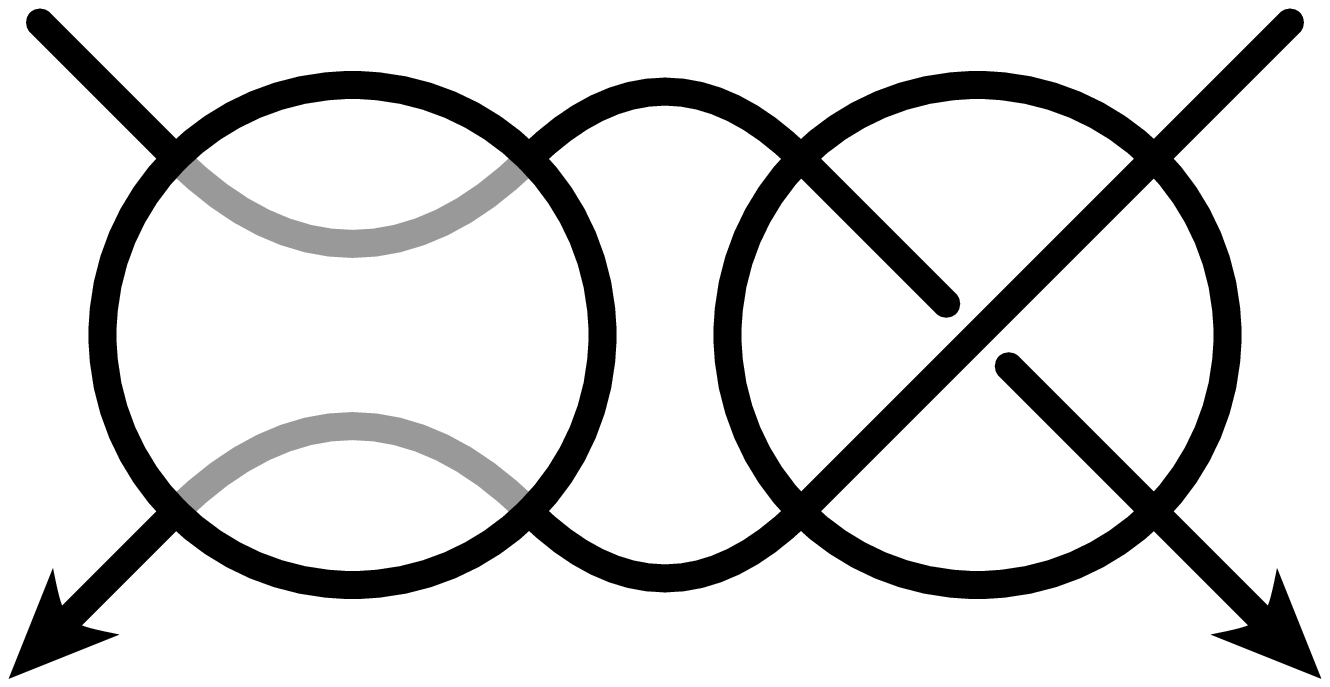}}\ar[r]^{c=-1} &
{\includegraphics[scale=0.1875]{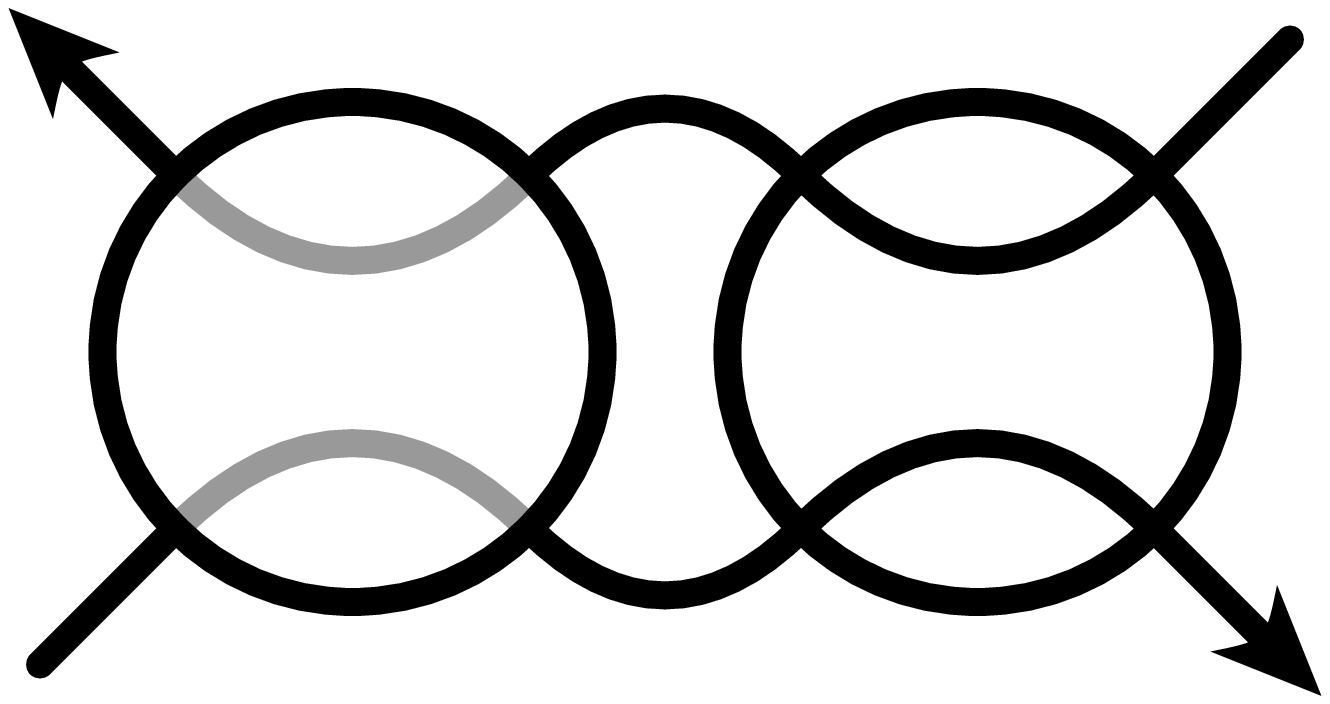}}
}\end{array}$}and a similar argument for the (switched) values of $c$. This time, we first resolve the (positive) crossing of $U^{\bar{\si}}$ to obtain $U$ with it's orientation unchanged. The induced orientation on $T^\si$ reverses (as before) the orientation of both strands so that the number of positive and negative crossings are once more unchanged. Therefore, we lose only a positive crossing, and obtain $c=n_-(U)-n_-(U^{\bar{\si}})=n_-\left(K_\beta(\positive\zero)\right)-n_-\left(K_\beta(\positive\negative)\right)=0$. When the distinguished crossing of $T^\si$ is resolved to obtain $T$, we have that the orientation on $U$ is once more preserved, so that the contribution to $c$ comes from comparing $T$ and $T^\si$ only. Once again, we remove the distinguished crossing (a negative crossing) and reverse orientation of both the strands of $T$. Therefore, $c=n_-(T)-n_-(T^\si)=n_-\left(K_\beta(\zero\zero)\right)-n_-\left(K_\beta(\positive\zero)\right)=-1$.

With this information we construct the second diagram of exact sequences, as in section \ref{Particular Example}:
\[\xymatrix@R=7pt@C=30pt{
	{\cdots\ L_{q-1}^{u-1}}\ar@/^1pc/[dr] &  & & &\\
{\cdots\ L_{q-1}^{u-1}}\ar[r] & {Kh^{u-1}_{q-2}K(\positive\zero)}\ar[r]\ar@/_1pc/[dr] & {Kh^u_q K(\positive\negative)}\ar[r] & {L^{u}_{q-1}\cdots\ } & {(2)}\\
& & {Kh^u_qK(\zero\zero )}\ar[r] & {L^{u}_{q-1}\cdots\ } &
}\]

We pause here to remark that the cases with different connectivity for $T^\si$ and $U^{\bar{\si}}$ ($\Gzero,\Gzero$ and $\Gzero,\Gcross$ and $\Gcross,\Gzero$) proceed in the same way, with only minor adjustments to the induced orientations. In fact, the proof amounts to reordering and/or rotating the oriented diagrams used above. We leave this step to the reader. Moreover, although the choice of $\beta = \si_1^{-1}\si_2\si_1^{-2n}$ will be sufficient for our purposes (cf. section \ref{Examples}) it is possible to consider other braids (with different associated permutations) by studying different cases of orientations induced by other permutations on the strands. Note also that acting by $\bar{\si}$ (instead of $\si$) on $K$ switches the two exact sequences in each of the diagrams (1) and (2).

With the diagrams (1) and (2) in hand the conclusion of the proof of lemma \ref{lem:khovanov} proceeds exactly as in section \ref{Particular Example}. In particular, it follows immediately that   
\[Kh_q^uK_\beta(T^\si,U^{\bar{\si}}) = Kh_q^uK_\beta(\positive\negative) \cong Kh^u_qK_\beta(\zero\zero) = Kh_q^uK_\beta(T,U)\]
whenever $u\ne0$ since $L^u=0$ whenever $u\ne0$. Applying lemma \ref{lem:jones} we have that the knots considered have the same (unormalized) Jones polynomial, which allows us to conclude that the free part is isomorphic for $u=0$, and applying the diagram (2) we conclude that the torsion part is isomorphic for $u=0$ as well since $L\cong(\bZ)_{-2}\pls(\bZ\pls\bZ)_0\pls(\bZ)_2$ is torsion free. In particular,
\[KhK_\beta(T^\si,U^{\bar{\si}})\cong KhK_\beta(T,U).\] 

\section{Examples}\label{Examples}

We have yet to see that the knots of lemma \ref{lem:khovanov} are distinct. There are three degrees of freedom in the construction: The choice of a braid 
$\be$ (that is, a choice of positive integer $n$), the choice of a (simple) tangle $T$, and the choice of a compatible (simple) tangle $U$ relative to some fixed $T$. This gives rise to a wide range of knots; we treat some particular examples. The notation for knots used below is consistent with Rolfsen's notation \cite{Rolfsen1976} for knots with fewer that 11 crossings, and \url{KNOTSCAPE} notation \cite{Thistlethwaite} otherwise (see also \cite{knotatlas}) up to mirrors. 

\subsection{Kanenobu's knots.}First set $n=1$ and consider the braid $\be=\si^{-1}_1\si_2\si^{-2}_1$. It is easy to check that a tangle composed of some number of horizontal twists (including none) is simple. This gives rise to a particular infinite family of knots with identical Khovanov homology. Let $K=K_\be\big(\Czero, \Czero\big)$ so that  $K^\si=K_\beta\big(\Czero^\si, \Czero^{\bar{\si}}\big)=K_\beta\big(\Cpositive, \Cnegative\big)$. The knots $K= 4_1\#4_1$, $K^\si= 8_9\ (=8a16)$ and $K^{\si^2}=12n462$ are shown below.
\begin{center}
\includegraphics[scale=0.25]{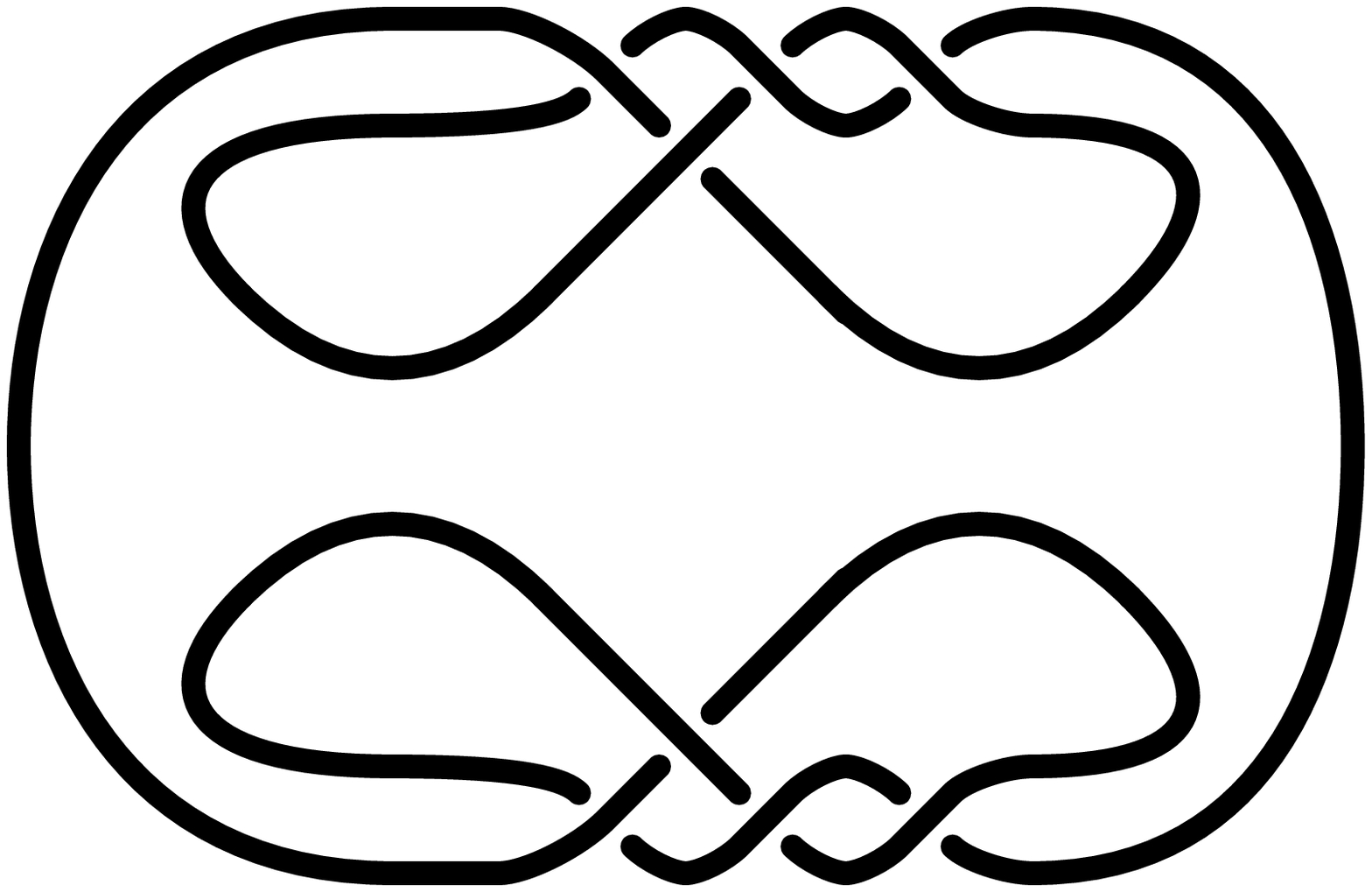}
\quad
\includegraphics[scale=0.25]{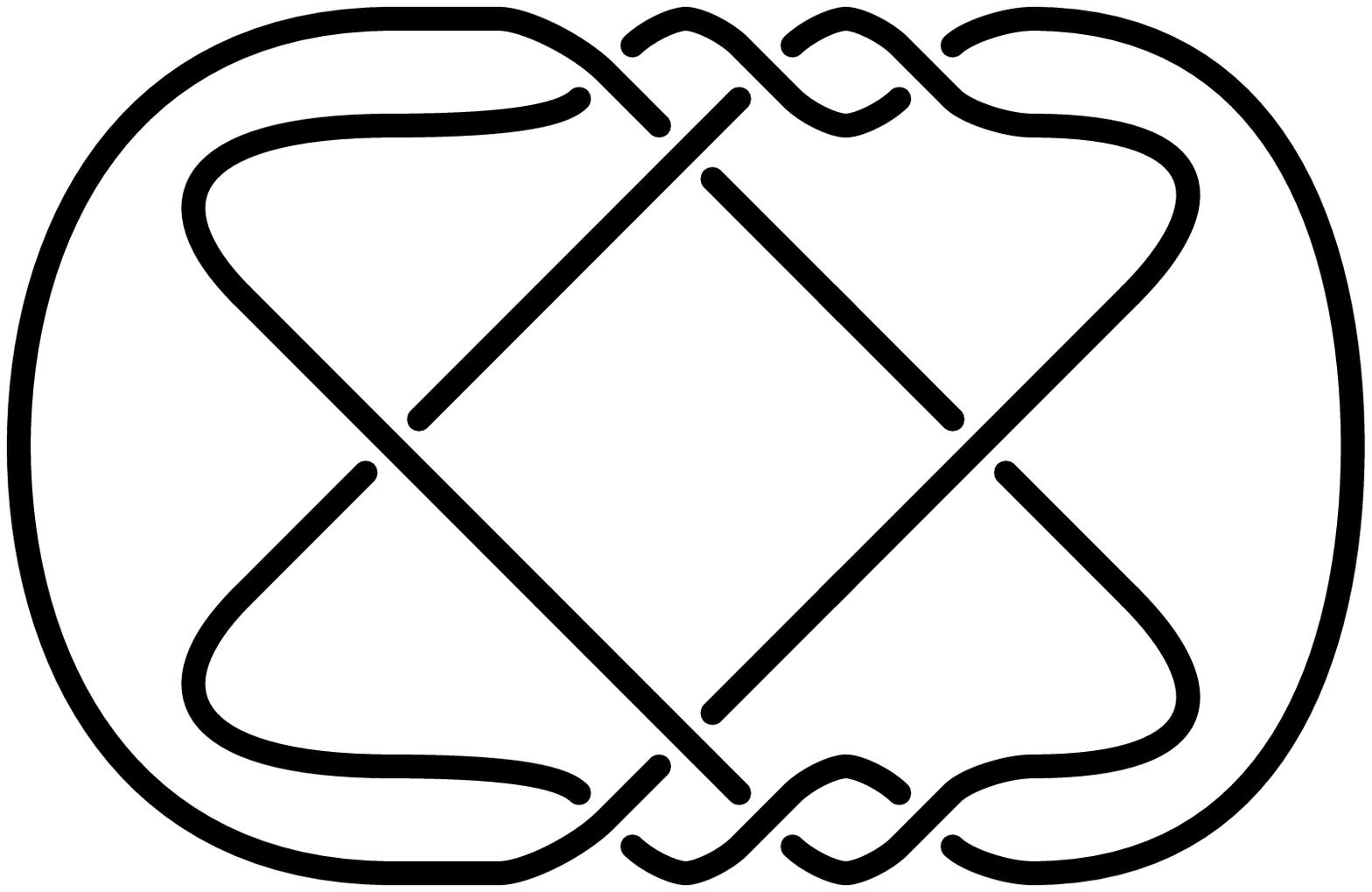}
\quad
\includegraphics[scale=0.25]{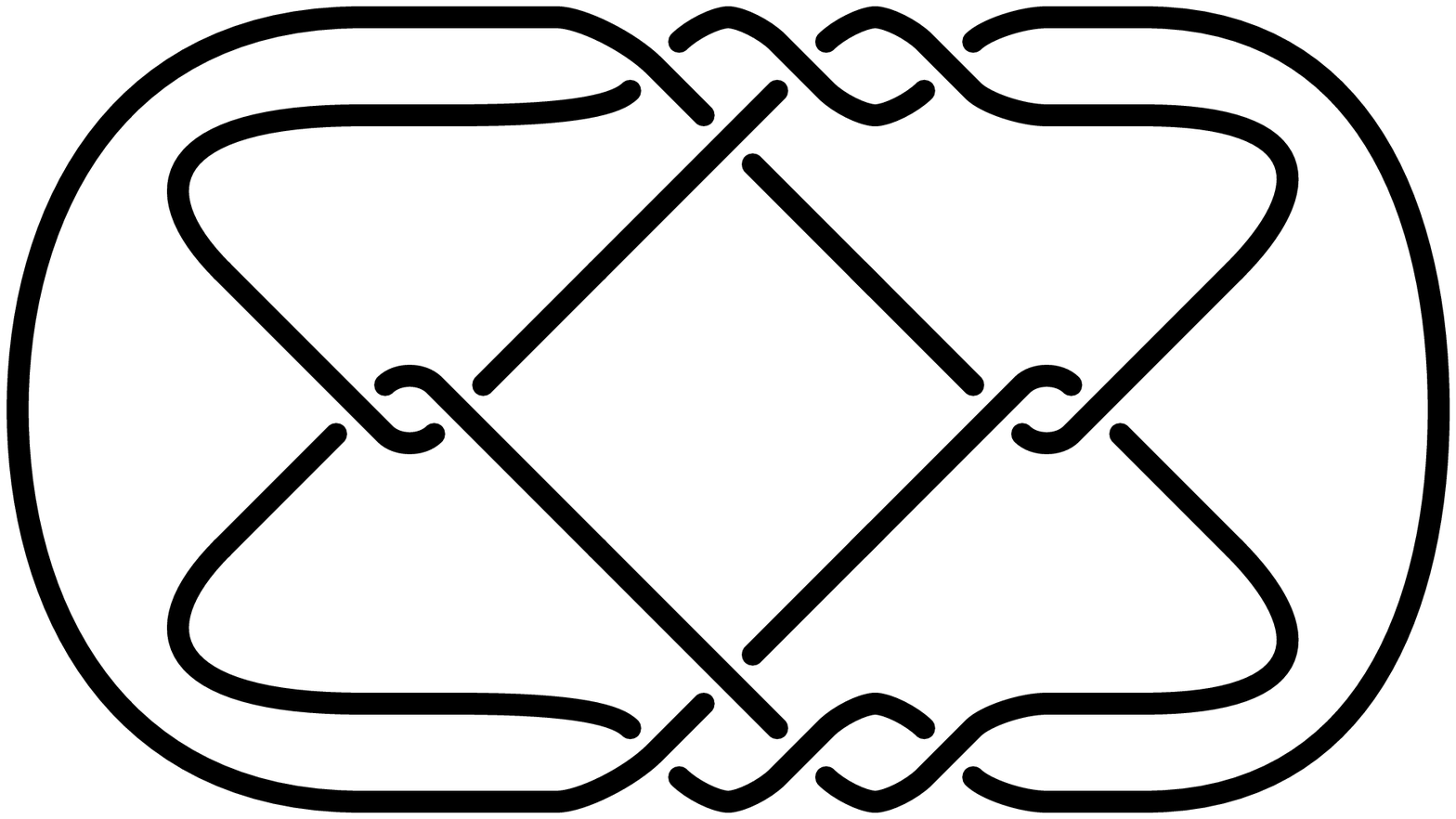}
\end{center}
By construction these knots have the same Khovanov homology, while $4_1\#4_1$ and $8_9$ (equivalently $8_9$ and $12n462$) have different HOMFLYPT polynomials. It should be noted that $4_1\#4_1$ and $12n462$ share the same HOMFLYPT polynomial, and the interested reader should consult \cite{Kanenobu1986} in which Kanenobu originally classified this example. In particular, this provides an infinite family of distinct knots with homology $Kh(8_9)$.    

\subsection{Distinct, prime knots with identical Khovanov homology.}
More generally, still using $\be=\si^{-1}_1\si_2\si^{-2}_1$, one can choose an arbitrary tangle $T$ and choose the mirror image $U=T^\star$ (which are always compatible) to obtain a pair of knots $K_\be(T,T^\star)$ and $K_\be(T^{\bar{\si}},(T^\star)^\si)$ that have the same Jones polynomial but distinct HOMFLYPT polynomial \cite{Watson2005}. As a result, they cannot be related by mutation \cite{Rolfsen1994}. Moreover, it can be shown that if $T$ is prime then each of the resulting knots are prime \cite{Watson2005}. A tangle $T=(B^3_T,\tau)$ is prime if and only if the two-fold branched cover of $B^3_T$ (branched over $\tau$) is irreducible and boundary irreducible \cite{Lickorish1981}. 
In the case where $T$ is also simple, we have proved the following:
\begin{theorem}\label{thm:pairs}
For every simple, prime tangle $T$ there exists a pair of distinct prime knots (each containing $T$) with identical Khovanov homology but distinct HOMFLYPT polynomial (and hence distinct triply-graded link homology).
\end{theorem}
\parpic[r]{$\begin{array}{c}\includegraphics[scale=0.25]{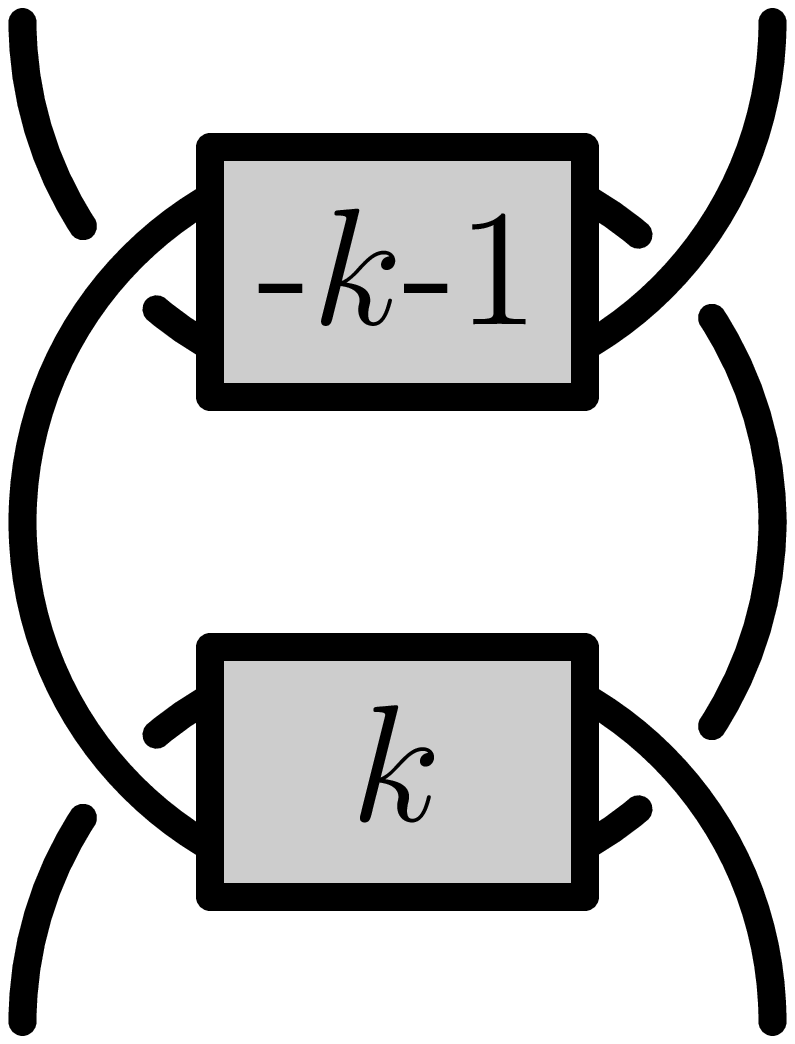}\end{array}$}
For example, the tangle on the right is prime for $k\ge0$ horizontal twists (cf. \cite{Lickorish1981}). It is also a simple tangle: $T+\Cone = \Cone$ by first applying $k+2$ Reidemeister type 2 moves, and then a single Reidemeister type 1 move. Thus, combining this tangle with theorem \ref{thm:pairs}, we get infinitely many pairs of prime knots that cannot be distinguished by Khovanov homology. Show below are the knots obtained from this construction in the case $k=0$. 

\begin{center}\includegraphics[scale=0.25]{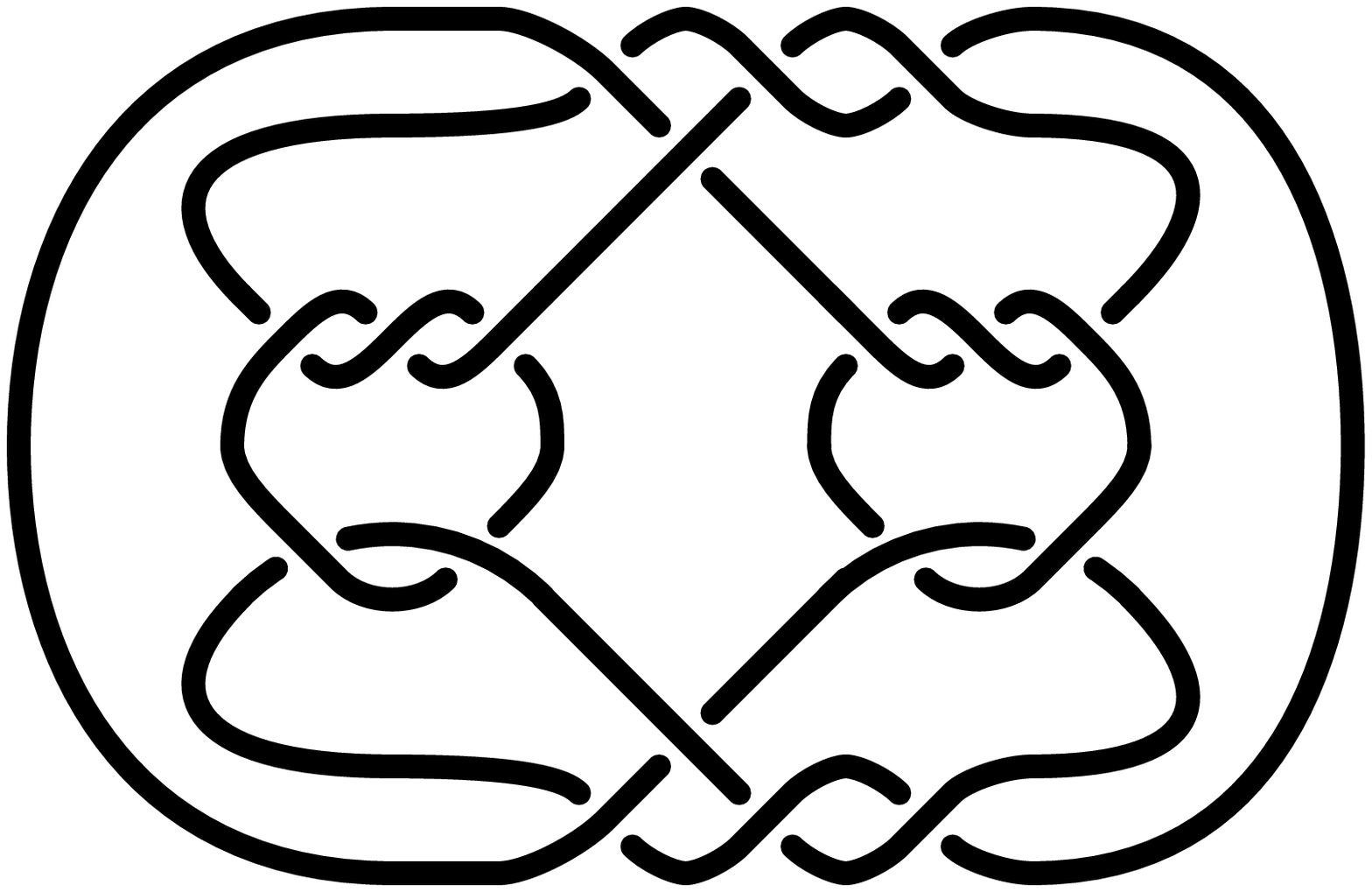}\quad\quad
\includegraphics[scale=0.25]{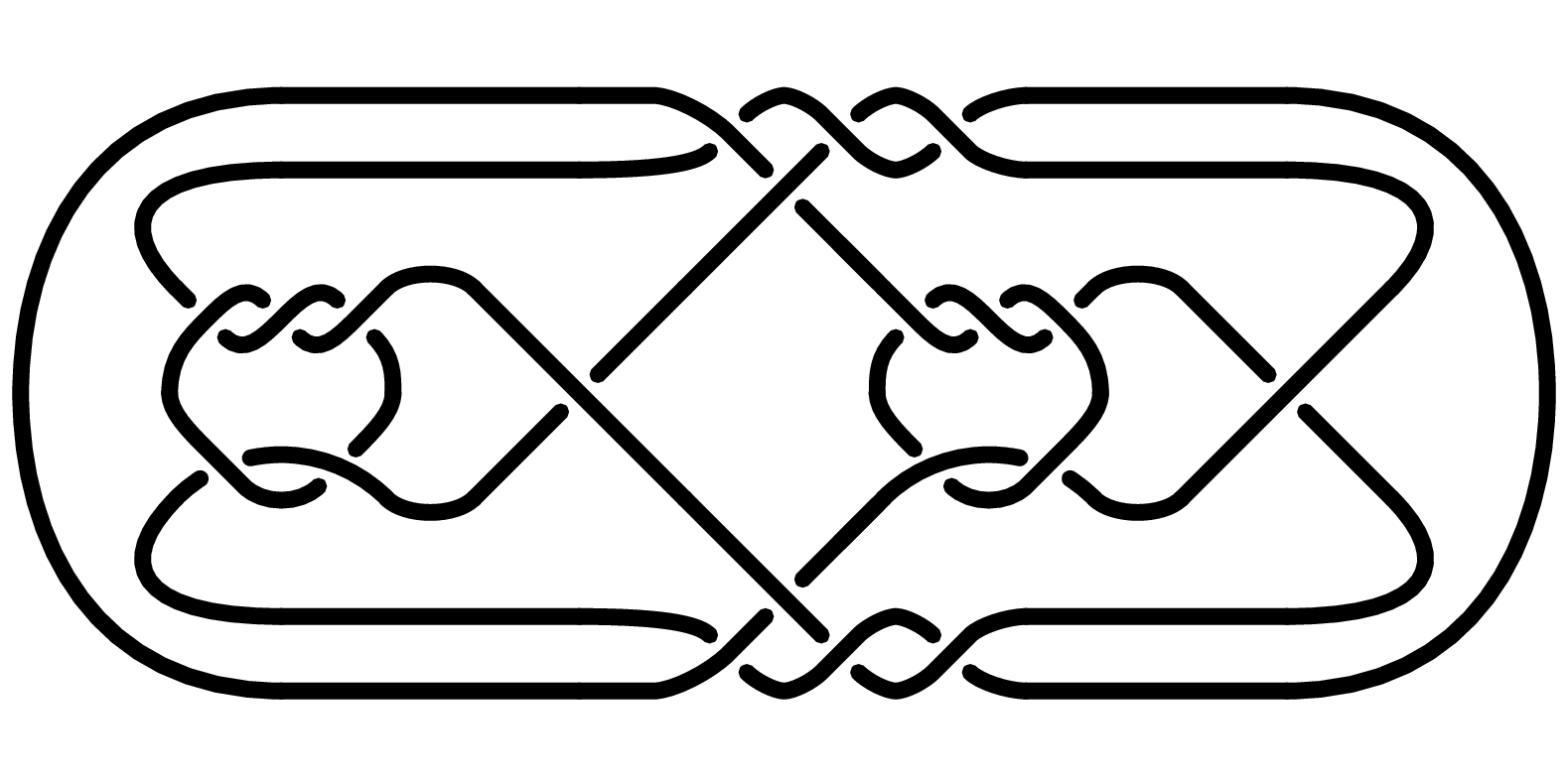}\end{center}

The fact that the pair of knots constructed in the proof of theorem \ref{thm:pairs} are distinguished by the HOMFLYPT depends on the fact that the basic pair $K= 4_1\#4_1$ and $K^\si= 8_9$ have distinct HOMFLYPT polynomial \cite{Watson2005}. Indeed, the proof goes through with any other pair of the form $K_\beta$ and $K_\beta^\si$ that are distinguished by the HOMFLYPT polynomial. Since, for example, $6_1\#6_1^\star$ and $12a1283$ (using $\beta=\si^{-1}_1\si_2\si^{-4}_1$) have distinct HOMFLYPT polynomials, we may revisit the construction in theorem \ref{thm:pairs} for this $\beta$ to obtain further examples of pairs of distinct prime knots that have identical Khovanov homology.

Note that since the pairs constructed by this method have different HOMFLYPT polynomial, they cannot be related by mutation.

\subsection{Constructing infinite families.} Although one needs a mechanism to prove that the knots obtained are distinct, the action of $\si$ defined in Section \ref{Construction} may be iterated (as in Kanenobu's examples) to obtain infinite families of knots with identical Khovanov homology. Luse and Rong classified the particular familly taking $\be=\si^{-1}_1\si_2\si^{-2n}_1$, and from this classification (cf. theorem 1.1 of \cite{LR2006}) we have the following:
\begin{theorem}\label{thm:infinite}
For each $n\in\bN$ there is an infinite family of distinct knots with identical Khovanov homology.
\end{theorem}
\begin{proof}
Fix $n\in\bN$ and consider the family of knots $K^l=K_\be(\si^{2l},\si^{-2l})$ where $\si^{2l}$ is the tangle consisting of $l$ horizontal full-twists (ie an element of the pure braid group on 2 strands), and $\be=\si^{-1}_1\si_2\si^{-2n}_1$ (note that in the case $n=1$ we have recovered Kanenobu's example \cite{Kanenobu1986}). As noted above, these tangles are simple, so by construction $K^l$ and $K^{l'}$ have identical Khovanov homology for any $l,l'\in\bZ$. Moreover, $K^l$ and $K^{l'}$ are distinct knots whenever $\gcd(l,2n+1)\ne\gcd(l',2n+1)$ \cite{LR2006}. If $p_1^{\al_1}p_2^{\al_2}\cdots p_k^{\al_k}$ is the prime decomposition of $2n+1$, we can choose $l=p_i$ (for any of the $i\in\{1,\ldots, k\}$) so that $\gcd(l,2n+1)=p_i$. Letting $l'$ range over all primes that do not appear in the prime decomposition of $2n+1$ gives the result. 
\end{proof}

\subsection{Odds and ends.}
As noted in the Section \ref{Proof}, further examples may be generated by altering the choice of $\be$. We give some more examples in this case. 

If $\beta=\si_1^{-1}\si_2\si_1^{-1}$ then $K=K_\be\big(\Czero,\Czero\big)=5_2\#5_2^\star$ and we obtain $K^\si=10_{48}=10a79$ and $K^{\si^2}=14n15498$ all of which have the same Khovanov homology, while $K$ and $K^\si$ are distinguished by the HOMFLYPT polynomial.

If $\beta=\si_1^{-2}\si_2\si_1^{-2}\si_2^{-1}\si_1$ then  $K=K_\be\big(\Czero,\Czero\big)=6_3\#6_3$ and we obtain $K^\si=12a819$ $K^{\si^2}=16n532490$ all of which which have the same Khovanov homology, while $K$ and $K^\si$ are distinguished by the HOMFLYPT polynomial.

While such $\beta$ were not treated here to streamline the proof of lemma \ref{lem:khovanov}, note that once again these examples generate pairs distinguished by the HOMFLYPT polynomial as in theorem \ref{thm:pairs}. 


Finally, we note that in many cases knots admitting a diagram of the form given in section \ref{Construction} have the same Khovanov homology even when the simplicity requirement on the tangles is dropped from the hypothesis of lemma \ref{lem:khovanov}. As an example, the knots $12a990$ and $12a1225$ arise in this way (consider the braid closure of $\be\si_2^3\bar{\be}\si_2^{-3}$ where $\beta=\si_1^{-1}\si_2\si_1^{-2}$). A second example of this phenomenon is given by the knots $12a427$ and $15n45009$ shown below. 
\begin{center}\includegraphics[scale=0.25]{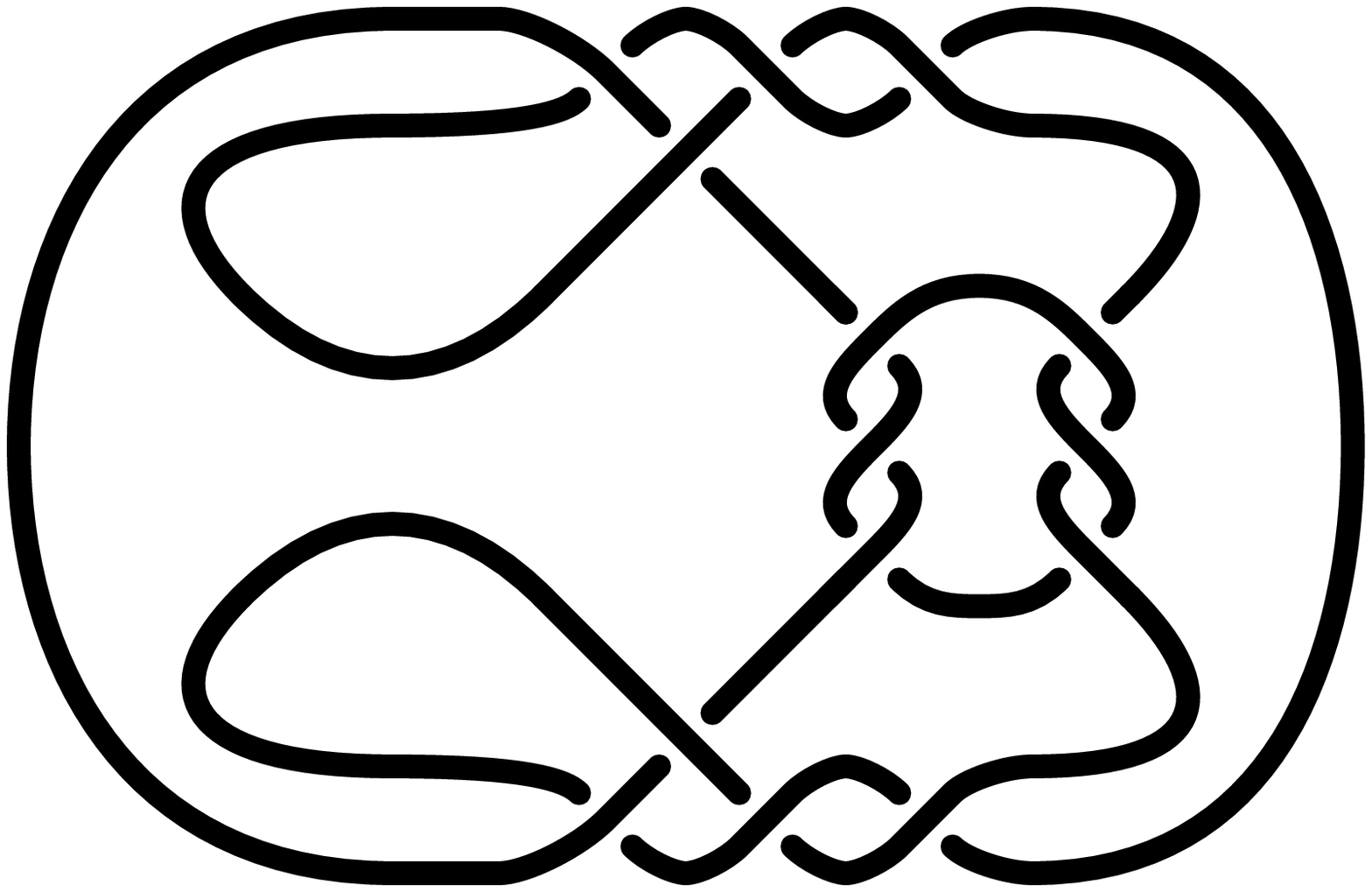}\quad\quad
\includegraphics[scale=0.25]{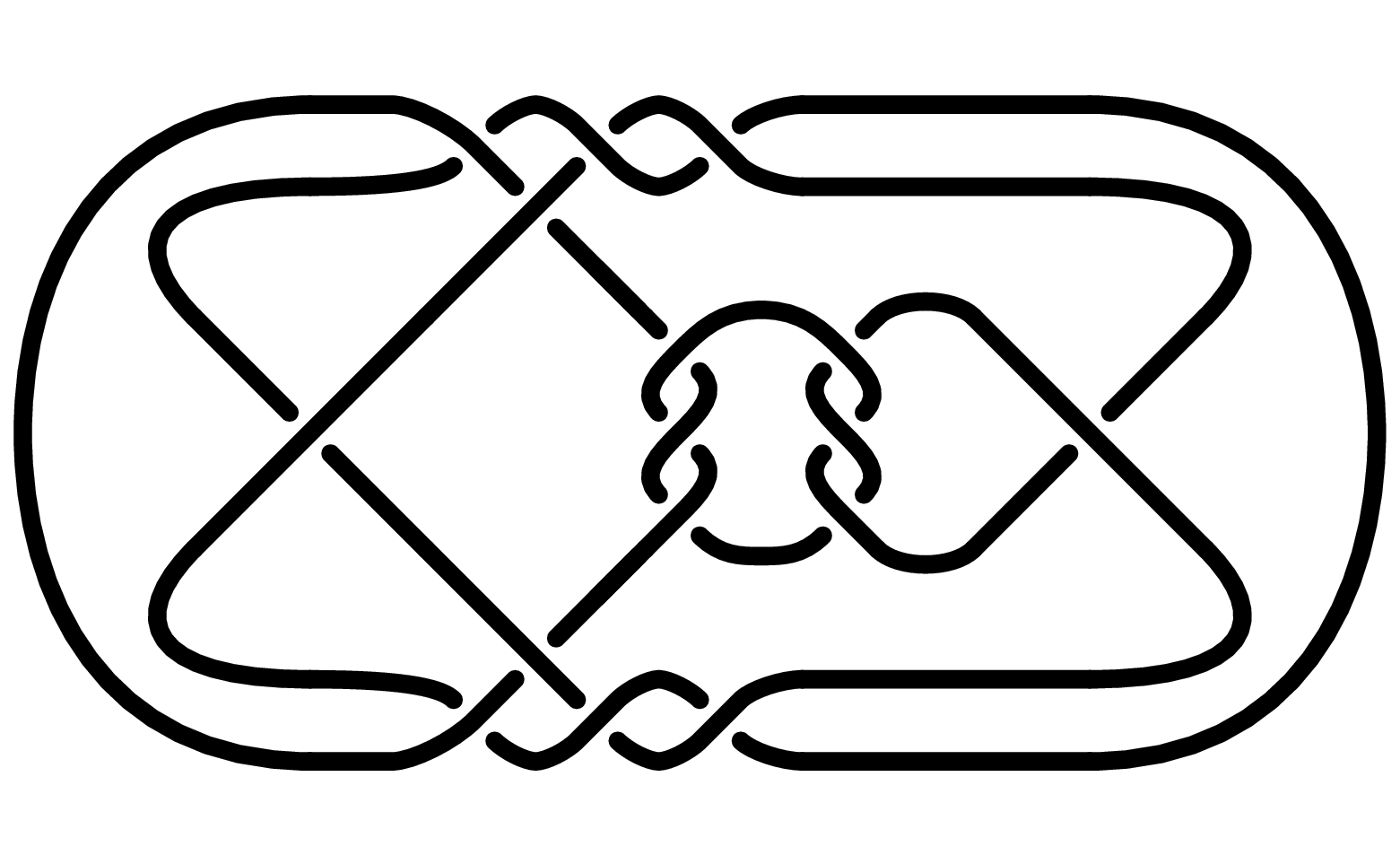}\end{center}
This suggests that there is another reason that such examples should have identical Khovanov homology, since the techniques used here cannot treat the most general case. 

\section{A Remark on Mutants}\label{Mutants}

\parpic[r]{$\begin{array}{c}\xymatrix@C25pt@R10pt{ 
{\includegraphics[scale=0.1875]{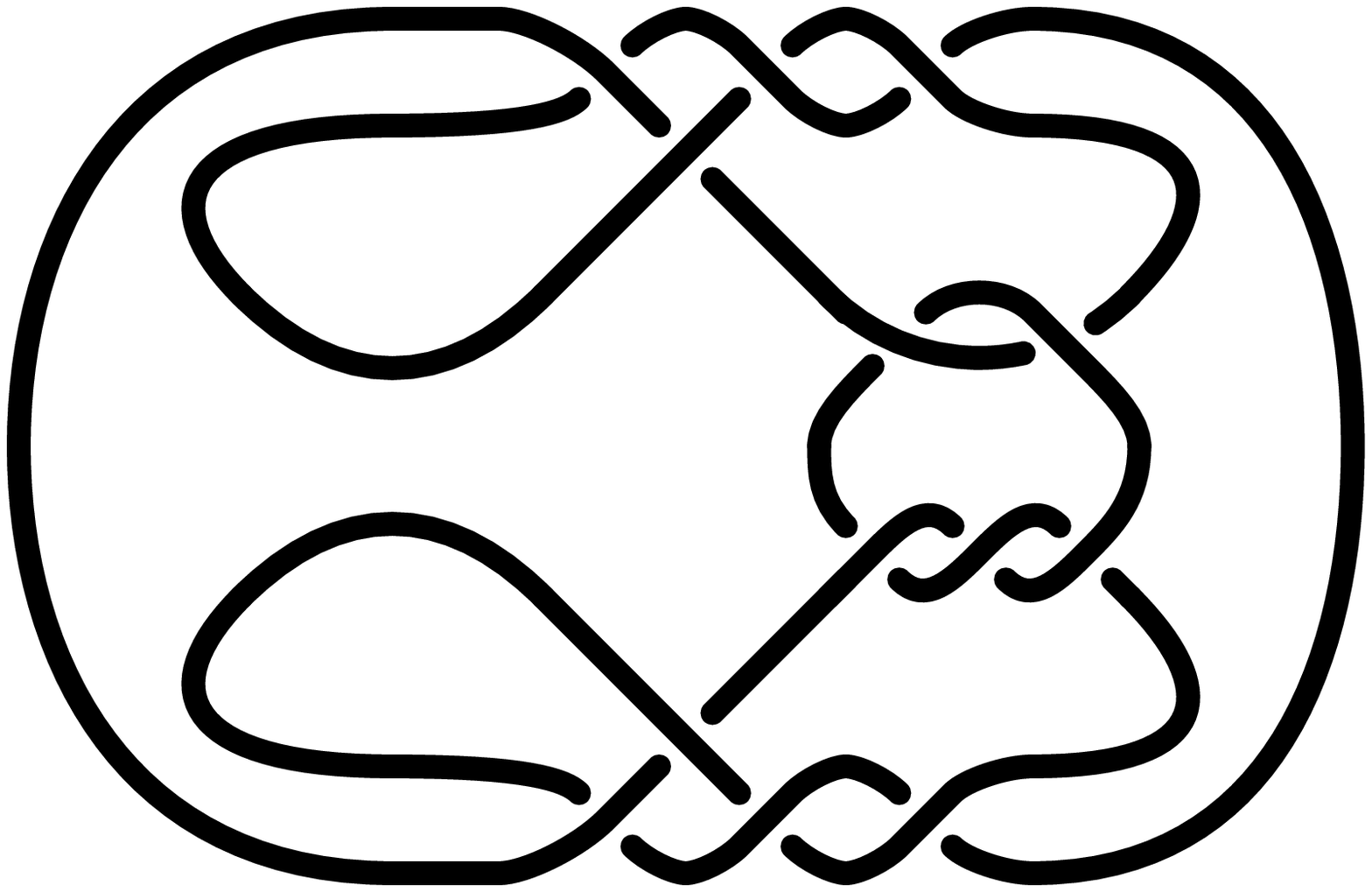}}\ar[r]^-{\bar{\si}}\ar[d]_-{{\rm mutate}} &
{\includegraphics[scale=0.1875]{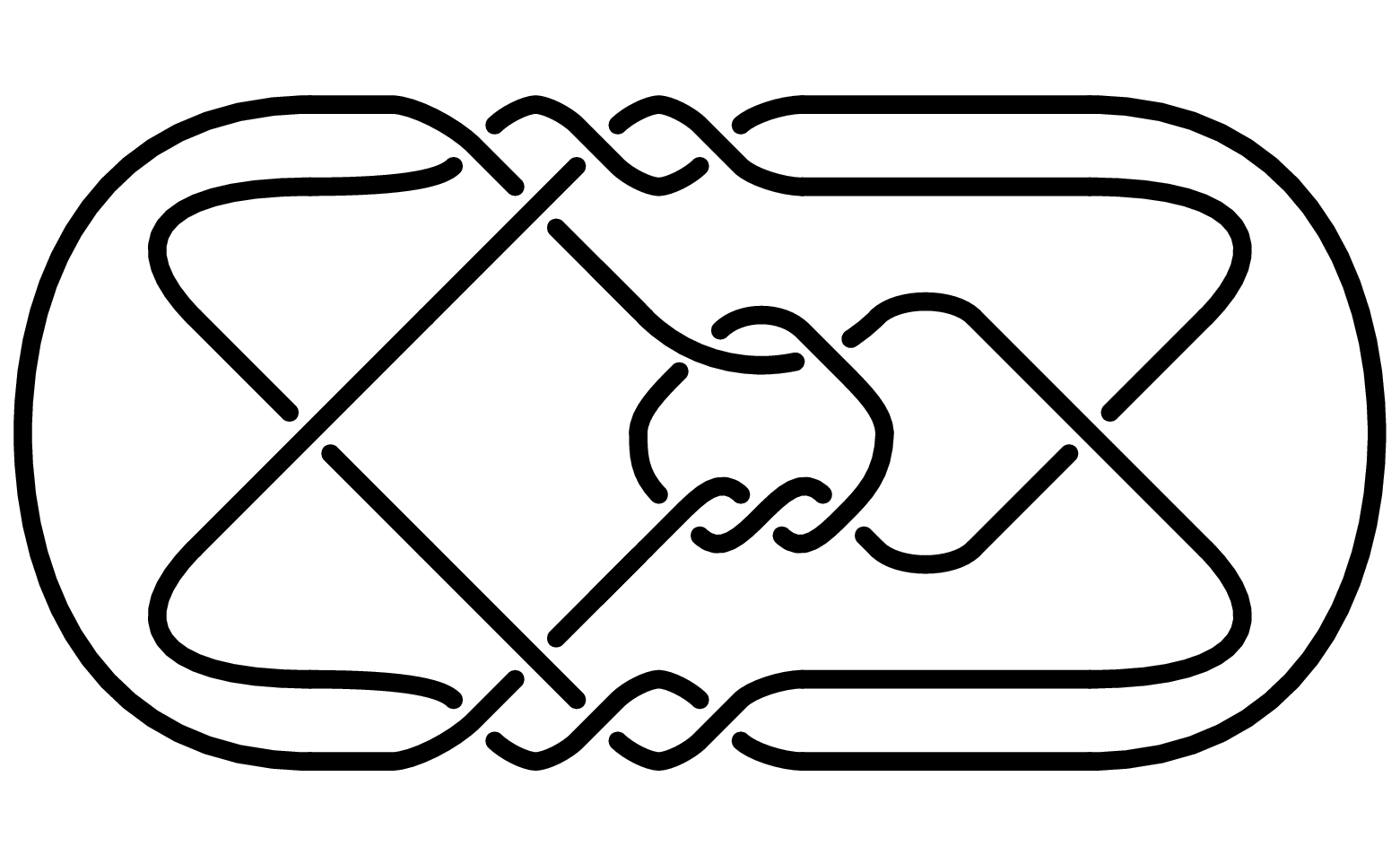}}\ar[d]^-{{\rm mutate}}\ar@{..>}[dl] \\
{\includegraphics[scale=0.1875]{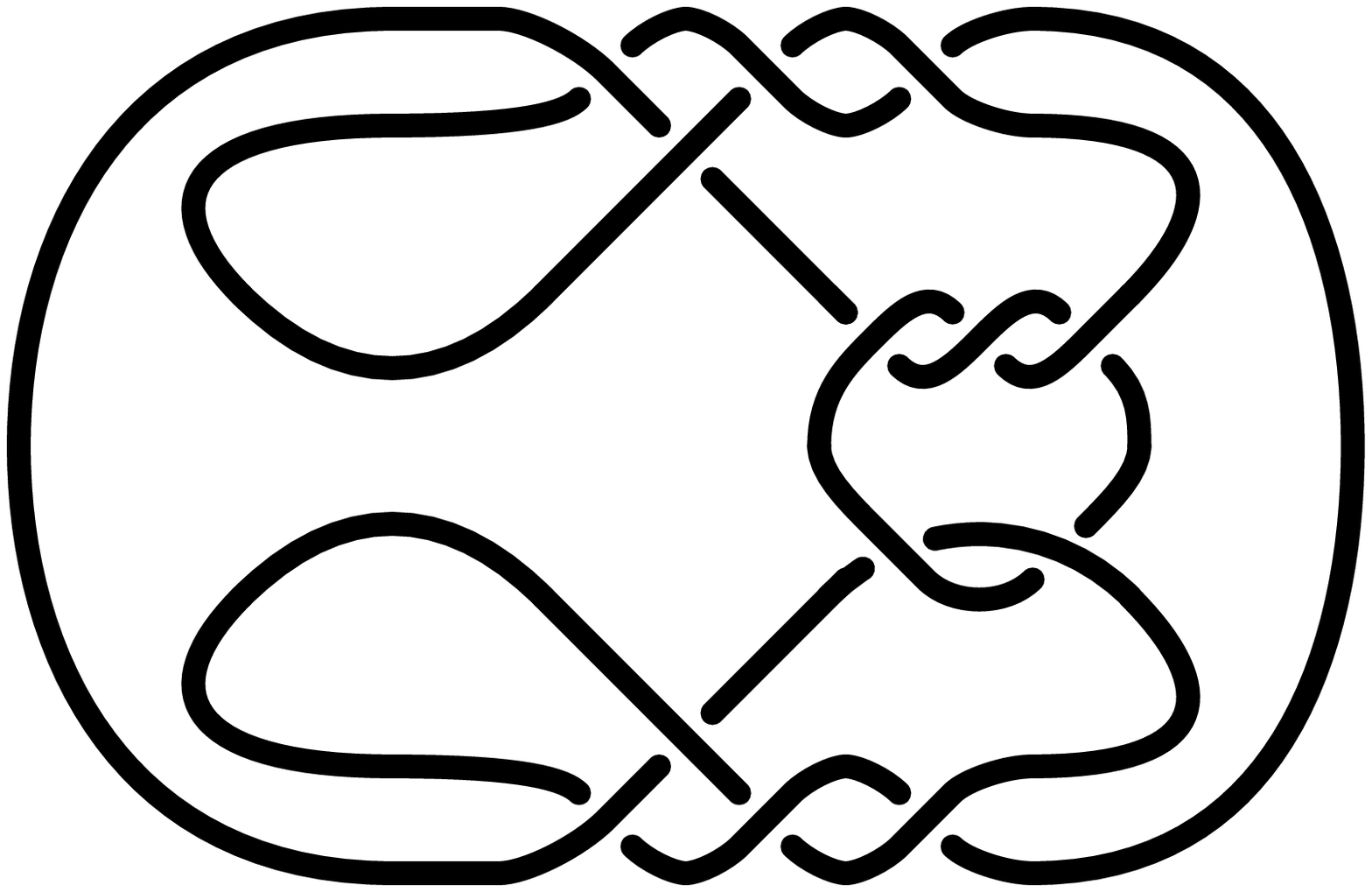}}\ar[r]^-{\bar{\si}} &
{\includegraphics[scale=0.1875]{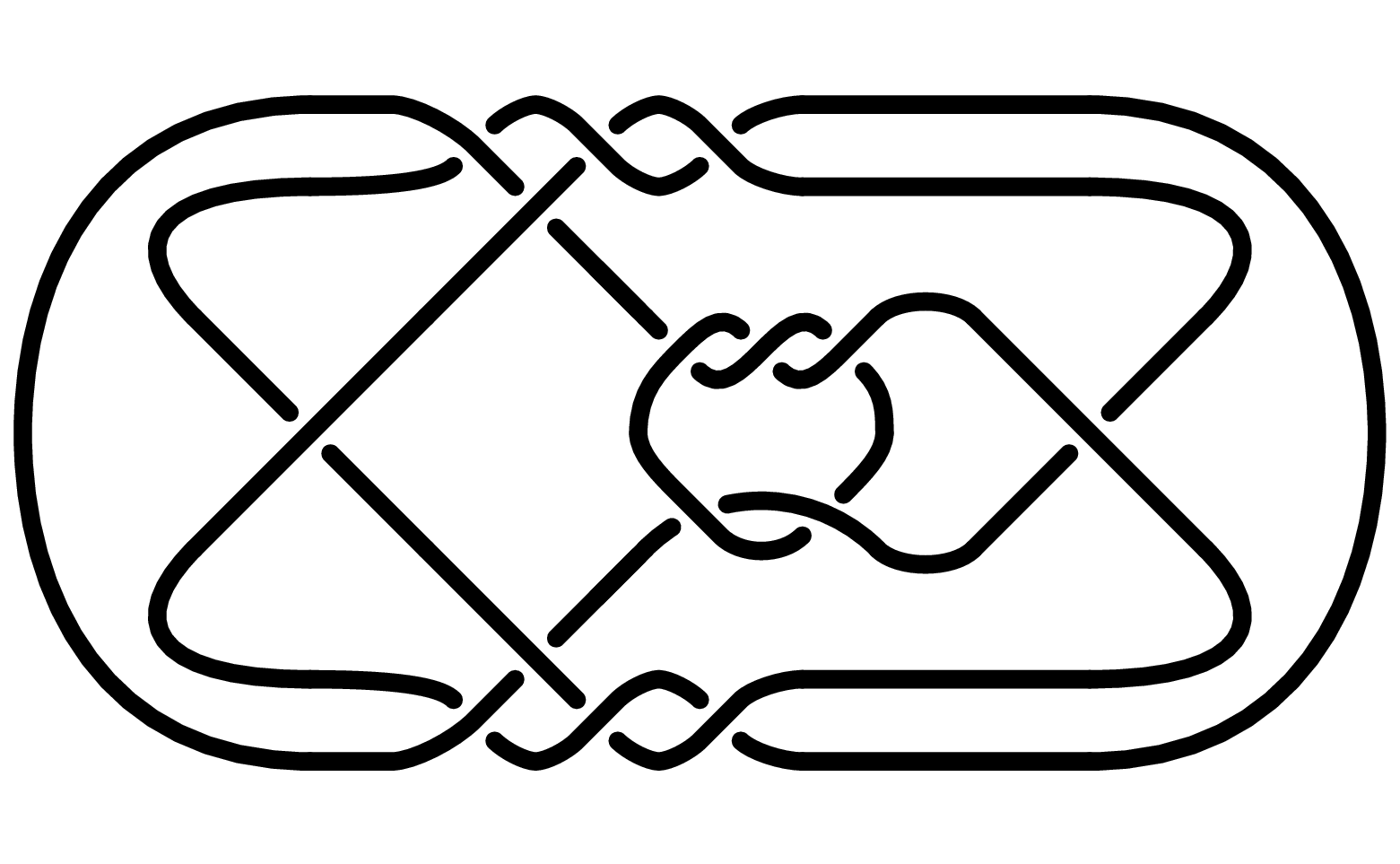}}
}\end{array}$}
We conclude with an interesting family of non-alternating knots with 13 crossings. The four knots in question may be arranged as in the diagram on the right, so that pair in the first row ($13n164$ and $13n922$) are related by twisting (and hence have identical Khovanov homology), as are the pair in the second row ($13n161$ and $13n795$). It may also be of interest to note that the (common) homology for all four of these knots is supported in 3 diagonals. The columns are related by mutation (see \cite{Rolfsen1994} for a general overview of mutation), and although it is unknown whether mutation preserves Khovanov homology for knots (cf. \cite{Bar-NatanWIKI,Wehrli2003}), it is possible to give an explanation for the phenomenon in this case. If we consider a similar construction to section \ref{Construction} allowing $B_2$ to act on the left (ie $\si: T \mapsto \sigleft$ and $\bar{\si}:T \mapsto \SIGleft$), the the diagonal arrow in the diagram corresponds to this action. That is, the mutant pair can be seen as the composition of right (twist) action, followed by a left (untwist) action. The proof of lemma \ref{lem:khovanov} goes through in the same way for this left action (on the same class of knots), and hence leaves the Khovanov homology invariant. 

In fact, we can make a slightly more general statement from this observation: If $K_\be\big(\Czero, T\big)$ (equivalently $K_\be\big(T,\Czero)$) then, the mutation that flips $T$ across the horizontal axis leaves the Khovanov homology invariant. We remark that this mutation is relative to the particular complimentary tangle (or external wiring \cite{Rolfsen1994}) $K_\be\big(\Czero, -\big)$, and requires, of course, that $T$ be simple. Despite these restrictions, it is interesting to note that this construction produces an infinite class of mutants that cannot be detected by Khovanov homology by altering either the braid $\beta$ or the (simple) tangle $T$.    

\bibliographystyle{plain}
\bibliography{paperV3}

\end{document}